\documentclass[11pt]{article}
\usepackage{amssymb,latexsym}
\usepackage{epsfig}
\usepackage{eufrak}
\usepackage{amsmath}
\usepackage{mathrsfs}
\usepackage{color}

\newtheorem{theorem}{Theorem}
\newtheorem{lemma}{Lemma}
\newtheorem{corollary}{Corollary}

\newcommand{\be}{\begin{equation}}
\newcommand{\ee}{\end{equation}}
\newcommand{\bee}{\begin{eqnarray*}}
\newcommand{\eee}{\end{eqnarray*}}
\newcommand{\bel}{\begin{eqnarray}}
\newcommand{\eel}{\end{eqnarray}}
\newcommand{\bec}{\begin{cases}}
\newcommand{\eec}{\end{cases}}
\newcommand{\bem}{\begin{bmatrix}}
\newcommand{\eem}{\end{bmatrix}}

\newcommand{\la}{\label}
\newcommand{\li}{\left}
\newcommand{\ri}{\right}

\newcommand{\DEF}{\stackrel{\mathrm{def}}{=}}

\newcommand{\ovl}{\overline}
\newcommand{\udl}{\underline}

\newcommand{\lc}{\lceil}
\newcommand{\rc}{\rceil}
\newcommand{\lf}{\lfloor}
\newcommand{\rf}{\rfloor}

\newcommand{\vep}{\varepsilon}

\newcommand{\de}{\delta}

\newcommand{\se}{\theta}
\newcommand{\Se}{\Theta}

\newcommand{\ze}{\zeta}
\newcommand{\al}{\alpha}
\newcommand{\ba}{\beta}

\newcommand{\ro}{\rho}

\newcommand{\om}{\omega}
\newcommand{\Om}{\Omega}

\newcommand{\f}{\frac}
\newcommand{\sq}{\sqrt}
\newcommand{\cd}{\cdots}

\newcommand{\qu}{\quad}
\newcommand{\qqu}{\qquad}
\newcommand{\fa}{\forall}

\newcommand{\mscr}{\mathscr}
\newcommand{\mcal}{\mathcal}
\newcommand{\mbf}{\mathbf}
\newcommand{\bb}{\mathbb}

\newcommand{\wh}{\widehat}

\newcommand{\mrm}{\mathrm}
\newcommand{\bs}{\boldsymbol}

\newcommand{\ap}{\approx}

\newcommand{\sh}{\slash}

\newcommand{\tx}{\text}

\newcommand{\iy}{\infty}

\newcommand{\bed}{\begin{description}}
\newcommand{\eed}{\end{description}}
\newcommand{\bei}{\begin{itemize}}
\newcommand{\eei}{\end{itemize}}
\newcommand{\ben}{\begin{enumerate}}
\newcommand{\een}{\end{enumerate}}
\newcommand{\bib}{\bibitem}
\newcommand{\beL}{\begin{lemma}}
\newcommand{\eeL}{\end{lemma}}
\newcommand{\beT}{\begin{theorem}}
\newcommand{\eeT}{\end{theorem}}
\newcommand{\beC}{\begin{corollary}}
\newcommand{\eeC}{\end{corollary}}
\newcommand{\sect}{\section}

\newcommand{\bpf}{\begin{pf}}
\newcommand{\epf}{\end{pf}}
\newcommand{\bsk}{\bigskip}

\newcommand{\pfbox}{\hfill\mbox{$\Box$}}

\newenvironment{pf}{\paragraph*{Proof{\rm.}}}{\pfbox\bigskip}

\begin{document}

\title{{\bf Sequential Tests of Statistical Hypotheses with Confidence Limits}
\thanks{The author had been previously working with
Louisiana State University at Baton Rouge, LA 70803, USA, and is now with Department of Electrical Engineering, Southern University and A\&M
College, Baton Rouge, LA 70813, USA; Email: chenxinjia@gmail.com. The main results of this paper have been presented in Proceedings of SPIE
Conferences, Orlando, April 5-9, 2010 and April 25-29, 2011.  The statistical methodology proposed in this paper has been applied to electrical
engineering and computer science, see recent literature \cite{Chen_SPIE, Chen_SPIE11a, Chen_SPIE11b, Chen_SPIE11c, Chen_J} and the references
therein. } }

\author{Xinjia Chen}

\date{First submitted in July 2010}

\maketitle

\begin{abstract}

  In this paper, we propose a general method for testing composite hypotheses.
  Our idea is to use confidence limits to define stopping and decision rules.
  The requirements of operating characteristic function can be
  satisfied by adjusting the coefficients of the confidence limits.
  For common distributions, such adjustment can be done via efficient
  computation by making use of the monotonicity of the associated operating characteristic
  function.  We show that the problem of testing multiple hypotheses can be cast into the general framework of constructing
  sequential random intervals with prescribed coverage probabilities.  We propose an inclusion principle for constructing multistage testing
  plans.  It is demonstrated that our proposed testing plans can be substantially more
  efficient than the sequential probability ratio test and its
  variations.  We apply our general methodology to develop an exact approach
   for testing hypotheses regarding the difference of two binomial proportions.

\end{abstract}

\section{Introduction}

A fundamental issue of statistical inference is to test the
parameter of a distribution.  Let $X$ be a random variable defined
in a probability space $(\Om, \mscr{F}, \Pr )$.  Assume that the
distribution of $X$ is determined by a parameter $\se \in \Se$,
where $\Se$ is the parameter space.  A frequent problem is to test
hypotheses $\mscr{H}_0: \se \leq \se_0$ versus $\mscr{H}_1: \se \geq
\se_1$ based on observations $X_1, X_2, \cd$ of $X$, where $\se_0 <
\se_1$ are two parametric values in $\Se$.  To control the risk of
making wrong decisions, it is typically required that \bel  &   &
\Pr \{ \tx{Reject} \; \mscr{H}_0 \mid \se \} \leq \al, \qqu \fa \se
\in \Se \; \tx{such that} \; \se \leq \se_0 \la{req1}\\
&   & \Pr \{ \tx{Accept} \; \mscr{H}_0 \mid \se \} \leq \ba, \qqu
\fa \se \in \Se \; \tx{such that} \; \se \geq \se_1, \la{req2} \eel
where $\al, \ba \in (0, 1)$. The probability $\Pr \{ \tx{Accept} \;
\mscr{H}_0 \mid \se \}$ for $\se \in \Se$ is referred to as the
Operating Characteristic (OC) function.  Since there is no
requirement imposed on $(\se_0, \se_1)$ for the OC function, such
interval $(\se_0, \se_1)$ is referred to as an {\it indifference
zone}.

A famous solution to the above problem is the Sequential Probability
Ratio Test (SPRT) proposed by Abraham Wald \cite{Wald} during World
War II of last century. Various modifications of SPRT have also been
proposed (see, \cite{Ghosh} and the references therein).  A major
drawback of SPRT and its variations is that the test plans are
efficient only for some parametric values. Another unpleasant
feature is that the sample number of the tests may not be bounded.
In this paper, to overcome the limitations of SPRT and its
variations, we take a different path. Instead of using likelihood
ratio as existing testing plans, we use confidence limits to define
stopping rules and decision rules.  A key idea is to compare
confidence limits with the endpoint of indifference zone. We
discovered that by adjusting the coefficients of the confidence
limits, it is possible to control the risk of making erroneous
decisions.  We demonstrate that, for testing problems related to
common distributions, the adjustment can be done by efficient
computation and for such purpose it suffices to evaluate the OC
function for the endpoints of the indifference zone.

The remainder of this paper is organized as follows.  In Section 2, we proposed our general principle of hypothesis testing in a sequential
framework.  In Section 3, we discuss the applications of the general principle to test the expectation of a random variable. Our development is
motivated by the fact that the means of many important random variables are exactly the parameters of their distributions. In Section 4, we
derived multistage  testing plans from the fully sequential testing plans.   In Section 5, we propose to use confidence limits to design
multistage plans for testing multiple hypotheses.  In Section 6, we apply the design method based on confidence limits to test hypotheses
regarding functions of two binomial proportions.  Section 7 is the conclusion. All proofs of theorems are given in Appendices.

Throughout this paper, we shall use the following notations.  The notation $\emptyset$ denotes an empty set. The set of real numbers is denoted
by $\bb{R}$.  The set of positive integers is denoted by $\bb{N}$.  The expectation of a random variable is denoted by $\bb{E}[.]$. We use the
notation $\Pr \{ . \mid \se \}$ to denote the probability of an event which is defined in terms of random variables parameterized by $\se$. The
parameter $\se$ in $\Pr \{ . \mid \se \}$  may be dropped whenever this can be done without introducing confusion.  If $Z$ is parameterized by
$\se$, we denote $\Pr \{  Z \leq z \mid \se \}$ by $F_Z(z, \se)$ and $\Pr \{ Z \geq z \mid \se \}$ by $G_Z(z, \se)$ respectively.    The support
of a random variable $Z$ is denoted by $I_Z$, i.e., $I_Z = \{ Z(\om): \om \in \Om \}$.  The standard normal distribution is denoted by
$\Phi(.)$.  The other notations will be made clear as we proceed.

\section{Fully Sequential Tests of One-sided Hypotheses}

In this section, we shall consider the problem of testing hypotheses
$\mscr{H}_0: \se \leq \se_0$ versus $\mscr{H}_1: \se \geq \se_1$
stated in the introduction.

 Let $\mscr{L} (X_1, \cd, X_n; \de)$ and $\mscr{U} (X_1, \cd, X_n; \de)$
 be confidence limits such that
\be \la{CI0}
 \Pr \{ \mscr{L} (X_1, \cd, X_n; \de)  < \se \mid \se \}
\geq 1 - \de, \qqu \Pr \{ \mscr{U} (X_1, \cd, X_n; \de) > \se \mid
\se \} \geq 1 - \de \ee for $\de \in (0, 1)$ and $n = 1, 2, \cd$.
Let $\ze$ be a positive
 number less than $\min (\f{1}{\al}, \f{1}{\ba})$.
 The stopping rule of our testing plan is defined as follows:

Continue sampling until $\mscr{L} (X_1, \cd, X_n; \ze \al) \geq
\se_0$ or $\mscr{U} (X_1, \cd, X_n; \ze \ba) \leq \se_1$ for some
$n$.

Upon termination of sampling, the decision rule of our testing plan
is defined as follows: \bed

\item [(i):] If $\se_0 \leq \mscr{L} (X_1, \cd, X_n; \ze \al) < \mscr{U}
(X_1, \cd, X_n; \ze \ba) \leq \se_1$ is not satisfied, then accept
$\mscr{H}_0$ when $\mscr{U} (X_1, \cd, X_n; \ze \ba) \leq \se_1$ and
reject $\mscr{H}_0$ when $\mscr{L} (X_1, \cd, X_n; \ze \al) \geq
\se_0$.

\item [(ii):] If $\se_0 \leq \mscr{L} (X_1, \cd, X_n; \ze \al) < \mscr{U}
(X_1, \cd, X_n; \ze \ba) \leq \se_1$ is satisfied, then accept
$\mscr{H}_0$ or reject $\mscr{H}_0$ based on an arbitrarily
pre-specified policy.  On possible way to specify such policy is to
accept $\mscr{H}_0$ if  \[ \f{ f (X_1, \cd, X_n; \se_0) }{ f (X_1,
\cd, X_n; \se_1) } \geq \f{\al}{\ba}
\]
and reject $\mscr{H}_0$ otherwise.

\eed

Here $f(x_1, \cd, x_n; \se)$ denotes the joint probability density
function (for the continuous case) or the joint probability mass
function (for the discrete case) of $X_1, \cd, X_n$ parameterized by
$\se$.

Throughout the paper, let $\mbf{n}$ denote the number of samples at
the termination of sampling. Clearly, if $\mscr{U} (X_1, \cd, X_n;
\ze \ba) - \mscr{L} (X_1, \cd, X_n; \ze \al)$ converges to $0$ in
probability as $n \to \iy$, then $\Pr \{ \mbf{n} < \iy \mid \se \} =
1$.   As will be seen in the next section, for a wide class of
hypothesis testing problems, we can show that the sample sizes are
absolutely bounded. Moreover, we can show that the probabilities of
committing decision errors for the above testing plan can be
adjusted below any prescribed level by choosing a sufficiently small
value of $\ze > 0$.

\sect{Testing the Expectation of a Random Variable}

In many situations, the parameter to be tested is equal to the
expectation of the associated random variable. That is, $\se =
\bb{E} [ X ]$. As proposed in the last section, our testing plan is
defined in terms of confidence limits.  In order to obtain
confidence limits and establish the monotonicity of the OC function
for the associated testing plan, we need to introduce the concept of
{\it unimodal-likelihood estimator} (ULE).   For a random tuple
$X_1, \cd, X_{\mbf{r}}$ (of deterministic or random length
$\mbf{r}$) parameterized by $\se$, we say that the estimator
$\varphi (X_1, \cd, X_{\mbf{r}})$ is a ULE of $\se$ if $\varphi$ is
a multivariate function such that, for any observation $(x_1, \cd,
x_{r})$ of $(X_1, \cd, X_{\mbf{r}})$, the likelihood function is
non-decreasing with respect to $\se$ no greater than $\varphi (x_1,
\cd, x_{r})$ and is non-increasing with respect to $\se$ no less
than $\varphi (x_1, \cd, x_{r})$. For discrete random variables
$X_1, \cd, X_{r}$, the associated likelihood function is the joint
probability mass function $\Pr \{ X_i = x_i, \; i = 1, \cd, r \mid
\se \}$. For continuous random variables $X_1, \cd, X_{r}$, the
corresponding likelihood function is, $f_{X_1, \cd, X_r} (x_1, \cd,
x_r, \se)$, the joint probability density function of random
variable $X_1, \cd, X_{r}$. It should be noted that a ULE may not be
a maximum-likelihood estimator (MLE). On the other side, a MLE may
not be a ULE.

In this section, we shall focus on random variables such that their
sample means are ULEs.  Let $\ovl{X}_n = \f{\sum_{i = 1}^n X_i}{n}$,
where $X_1, X_2, \cd$ are i.i.d. samples of $X$ parameterized by
$\se = \bb{E} [ X]$. Define
\[
F_{\ovl{X}_n} (z, \se) = \Pr \li \{  \ovl{X}_n \leq z \mid \se  \ri
\}, \qqu G_{\ovl{X}_n} (z, \se) = \Pr \li \{  \ovl{X}_n \geq z \mid
\se  \ri \}.
\]
Suppose  that $\li \{ \se \in \Se: F_{\ovl{X}_n} (z, \se) \leq \de
\ri \}$ and $\li \{ \se \in \Se: G_{\ovl{X}_n} (z, \se) \leq \de \ri
\}$ are non-empty for $z \in I_{\ovl{X}_n}, \; n \in \bb{N}, \; \de
\in (0, 1)$. Define \be \la{CIA}
 L (z, \de) = \max \li \{ \se \in
\Se: G_{\ovl{X}_n} (z, \se) \leq \de \ri \}, \qqu U (z, \de) = \min
\li \{ \se \in \Se: F_{\ovl{X}_n} (z, \se) \leq \de \ri \}. \ee
Therefore, (\ref{CI0}) is satisfied if we define \be \la{CIB}
\mscr{L} (X_1, \cd, X_n; \de) = L(\ovl{X}_n, \de), \qqu \mscr{U}
(X_1, \cd, X_n; \de) = U(\ovl{X}_n, \de). \ee Suppose that the
moment generating function $\mscr{M} (\ro, \se) = \bb{E} [e^{\ro
X}]$ of $X$ exists for any $\ro \in \bb{R}$. Define $\mscr{C} (z,
\se) = \inf_{\ro \in \bb{R} } e^{- \ro z} \bb{E} [e^{\ro X}]$.  For
testing plans described in Section 2 with confidence limits defined
by (\ref{CIA}) and (\ref{CIB}), we have the following results.

\beT \la{thm1}

Suppose that the sample mean $\ovl{X}_n = \f{ \sum_{i=1}^n X_i }{n}$
is a ULE of $\se$. Then, \[ \mbf{n}  < \max \li \{  \f{ \ln (\ze
\al) } { \ln \mscr{C} ( \f{\se_0 + \se_1}{2}, \se_0 ) },  \; \f{ \ln
(\ze \ba) } { \ln \mscr{C} ( \f{\se_0 + \se_1}{2}, \se_1 ) } \ri \}.
\]
Moreover, $\Pr \{ \tx{Accept} \; \mscr{H}_0 \mid \se \}$ is
non-increasing with respect to $\se \in \Se$ such that $\se \notin
(\se_0, \se_1)$. Furthermore, both $\Pr \{ \tx{Reject} \; \mscr{H}_0
\mid \se_0 \}$ and $\Pr \{ \tx{Accept} \; \mscr{H}_0 \mid \se_1 \}$
tend to $0$ as $\ze \to 0$. \eeT

See Appendix \ref{thm1_app}  for a proof. As can be seen from
Theorem \ref{thm1}, the risk requirements (\ref{req1}) and
(\ref{req2}) can be satisfied by choosing a sufficiently small $\ze
> 0$. For this purpose, it suffices to ensure that (\ref{req1}) and
(\ref{req2}) hold for the endpoints of the indifference zone. If $X$
is a Bernoulli random variable, then the confidence limits are
exactly the same as Clopper-Pearson \cite{Clopper}.   If $X$ is a
Poisson random variable, then the confidence limits are exactly the
same as \cite{Garwood}.   It should be noted that in these
particular cases, recursive algorithms can be developed for the
evaluation of the OC functions and consequently facilitate the
determination of $\ze$ satisfying (\ref{req1}) and (\ref{req2}).

Note that confidence limits (\ref{CIA}) and (\ref{CIB}) are directly
derived from the exact tail probabilities of $\ovl{X}_n$.  In many
cases, simple bounds for the tail probabilities of $\ovl{X}_n$ are
available.  Hence, it is reasonable to hope that the stopping
boundary of the testing plans may be simpler if we use confidence
limits obtained from the bounds of the tail probabilities of
$\ovl{X}_n$.  Toward this goal, we have the following results
regarding confidence limits.

\beT \la{thm2} Suppose that the sample mean $\ovl{X}_n = \f{
\sum_{i=1}^n X_i }{n}$ is a ULE of $\se$. Suppose that $\{ \se \in
\Se:   [ \mscr{C} (z, \se) ]^n \leq \de, \; \se \leq z \}$ and $\{
\se \in \Se: [ \mscr{C} (z, \se) ]^n \leq \de, \; \se \geq z \}$ are
nonempty for $z \in I_{\ovl{X}_n}, \; n \in \bb{N}$ and $\de \in (0,
1)$. Define
\[
\mcal{L} (z, \de) = \max \li \{ \se \in \Se:   [ \mscr{C} (z, \se)
]^n \leq \de, \; \se \leq z \ri \}, \qqu \mcal{U} (z, \de) = \min
\li \{ \se \in \Se:   [ \mscr{C} (z, \se) ]^n \leq \de, \; \se \geq
z \ri \}.
\]
Then, {\small \[ \Pr \{ \mcal{L} ( \ovl{X}_n, \de ) < \se \} \geq 1
- \de, \qqu  \Pr \{ \mcal{U} ( \ovl{X}_n, \de ) > \se \} \geq 1 -
\de, \qqu \Pr \{ \mcal{L} ( \ovl{X}_n, \de ) < \se < \mcal{U} (
\ovl{X}_n, \de ) \} \geq 1 - 2 \de.
\]}
Moreover, both $\mcal{L} (z, \de)$ and $\mcal{U} (z, \de)$ are
non-decreasing with respect to $z \in I_{\ovl{X}_n}$. \eeT

See Appendix \ref{thm2_app} for a proof.  For testing plans
described in Section 2 with confidence limits defined by Theorem
\ref{thm2}, we have the following results.

\beT \la{thm3} Suppose that the sample mean $\ovl{X}_n = \f{
\sum_{i=1}^n X_i }{n}$ is a ULE of $\se$.
 Then,  \[ \mbf{n} < \max \li \{ \f{ \ln (\ze \al)
} { \ln \mscr{C} (  \f{\se_0 + \se_1}{2}, \se_0 ) }, \; \f{ \ln (\ze
\ba) } { \ln \mscr{C} (  \f{\se_0 + \se_1}{2}, \se_1 ) } \ri \}.
\]
Moreover,  $\Pr \{ \tx{Accept} \; \mscr{H}_0 \mid \se \}$ is
non-increasing with respect to $\se \in \Se$ such that $\se \notin
(\se_0, \se_1)$.  Furthermore, both $\Pr \{ \tx{Reject} \;
\mscr{H}_0 \mid \se_0 \}$ and $\Pr \{ \tx{Accept} \; \mscr{H}_0 \mid
\se_1 \}$ tend to $0$ as $\ze \to 0$.

\eeT

See Appendix \ref{thm3_app} for a proof.  Although our proposed
testing plans are derived from confidence limits, we can avoid the
computation of confidence limits by virtue of the monotonicity of
the tail probabilities (or their bounds) of $\ovl{X}_n$ with respect
to $\se$.   Under the assumption that the sample mean $\ovl{X}_n =
\f{ \sum_{i=1}^n X_i }{n}$ is a ULE of $\se$, we have the following
observations:

\bed

\item
 (i): $L (\ovl{X}_n, \ze \al) \geq \se_0 \Leftrightarrow G_{\ovl{X}_n} (\ovl{X}_n, \se_0) \leq \ze
 \al$ and $U (\ovl{X}_n, \ze \ba) \leq \se_1 \Leftrightarrow F_{\ovl{X}_n} (\ovl{X}_n, \se_1) \leq \ze
 \ba$.

\item (ii): $\mcal{L} (\ovl{X}_n, \ze \al) \geq \se_0 \Leftrightarrow \mscr{C}
(\ovl{X}_n, \se_0) \leq \ze
 \al, \; \ovl{X}_n \geq \se_0$ and
 $\mcal{U} (\ovl{X}_n, \ze \ba) \leq \se_1 \Leftrightarrow \mscr{C} (\ovl{X}_n, \se_1) \leq \ze
 \ba, \; \ovl{X}_n \leq \se_1$.
 \eed

Here the notation ``$\Leftrightarrow$'' means ``if and only if''.

\sect{Multistage Tests of One-sided Hypotheses}

In the preceding discussion, we have been focusing on fully
sequential tests. In many applications, it may not be practical to
take one sample at a time. It may be much more efficient to perform
testing  with samples in groups.  This motivates us to adapt the
sequential procedure proposed in Section 2 to multistage tests.

Consider hypotheses $\mscr{H}_0: \se \leq \se_0$ versus $\mscr{H}_1:
\se \geq \se_1$ as before.  We shall test such hypotheses with a
testing plan of $s$ stages.  Let $\ze$ be a positive
 number less than $\min (\f{1}{\al}, \f{1}{\ba})$.
 Let $n_1 < n_2 < \cd < n_s$ be sample sizes, where $n_\ell$ is the
sample size of the $\ell$-th stage.  Assume that the largest sample
size $n_s$ is no less than  the smallest integer $n$ such that $\{
\mscr{L} (X_1, \cd, X_n; \ze \al) < \se_0 < \se_1 <
 \mscr{U} (X_1, \cd, X_n; \ze \ba) \} = \emptyset$.

The stopping rule of our testing plan is defined as follows:

Continue sampling until $\mscr{L} (X_1, \cd, X_{n_\ell}; \ze \al)
\geq \se_0$ or $\mscr{U} (X_1, \cd, X_{n_\ell}; \ze \ba) \leq \se_1$
for some $\ell \in \{1, \cd, s \}$.

Upon termination of sampling, the decision rule of our testing plan
is defined as follows:

\bed

\item [(i):] If $\se_0 \leq \mscr{L} (X_1, \cd, X_{n_\ell}; \ze \al) <
\mscr{U} (X_1, \cd, X_{n_\ell}; \ze \ba) \leq \se_1$ is not
satisfied, then accept $\mscr{H}_0$ when $\mscr{U} (X_1, \cd,
X_{n_\ell}; \ze \ba) \leq \se_1$ and reject $\mscr{H}_0$ when
$\mscr{L} (X_1, \cd, X_{n_\ell}; \ze \al) \geq \se_0$.

\item [(ii):]  If $\se_0 \leq \mscr{L} (X_1, \cd, X_{n_\ell}; \ze \al) <
\mscr{U} (X_1, \cd, X_{n_\ell}; \ze \ba) \leq \se_1$ is satisfied,
then accept $\mscr{H}_0$ or reject $\mscr{H}_0$ based on an
arbitrarily specified policy.  On possible option of such policy is
to accept $\mscr{H}_0$ if
\[
\f{ f (X_1, \cd, X_{n_\ell}; \se_0) }{ f (X_1, \cd, X_{n_\ell};
\se_1) } \geq \f{\al}{\ba}
\]
and reject $\mscr{H}_0$ otherwise.

\eed

 Here $f(x_1, \cd, x_{n_\ell}; \se)$ denotes the joint
probability density function for the continuous case or the joint
probability mass function for the discrete case.  The confidence
limits $\mscr{L} (X_1, \cd, X_{n_\ell}, \ze \al)$ and $\mscr{U}
(X_1, \cd, X_{n_\ell}, \ze \ba)$ can be defined by (\ref{CIA}),
(\ref{CIB}) or Theorem \ref{thm2}.  With regard to such plans for
testing $\se = \bb{E} [ X]$, we have the following results.

\beT \la{thm4} Suppose that the sample mean $\ovl{X}_n = \f{
\sum_{i=1}^n X_i }{n}$ is a ULE of $\se$.
 Then, $\Pr \{ \tx{Accept} \; \mscr{H}_0 \mid \se \}$ is non-increasing
with respect to $\se \in \Se$ such that $\se \notin (\se_0, \se_1)$.
Moreover, $\Pr \{ \tx{Reject} \; \mscr{H}_0 \mid \se_0 \} \leq s \ze
\al$ and $\Pr \{ \tx{Accept} \; \mscr{H}_0 \mid \se_1 \} \leq s \ze
\ba$. Furthermore, both $\Pr \{ \tx{Reject} \; \mscr{H}_0 \mid \se_0
\}$ and $\Pr \{ \tx{Accept} \; \mscr{H}_0 \mid \se_1 \}$ tend to $0$
as $\ze \to 0$.

\eeT

Theorem \ref{thm4} can be shown by similar methods as that of
Theorems \ref{thm1} and \ref{thm3}.

\section{Multistage Tests of Multiple Hypotheses}

In this section, we shall propose a unified approach for constructing
multistage plans for testing multiple hypotheses by virtue of sequences of
fixed-sample-size confidence intervals.  The general problem to be considered is described as follows.

Let $X$ be a random variable defined in a probability space $(\Om,
\mscr{F}, \Pr )$.  Suppose the distribution of $X$ is determined by
an unknown parameter $\se$ in a parameter space ${\Se}$. In many
applications, it is desirable to infer the true value of $\se$ from
random samples $X_1, X_2, \cd$ of $X$. This topic can be formulated
as a general problem of testing $m$ mutually exclusive and
exhaustive composite hypotheses: \be \la{mainpr} \mscr{H}_0: \se \in
{\Se}_0, \qu \mscr{H}_1: \se \in {\Se}_1, \qu \ldots, \qu
\mscr{H}_{m - 1}: \se \in {\Se}_{m - 1}, \ee where $\Se_0 = \{ \se
\in \Se: \se \leq \se_1 \}, \; \Se_{m-1} = \{ \se \in \Se: \se
> \se_{m - 1} \}$ and $\Se_i = \{ \se \in \Se: \se_i < \se
\leq \se_{i+1} \}, \; i = 1, \cd, m - 2$ with $\se_1 < \se_2 < \cd <
\se_{m - 1}$.  To control the probabilities of making wrong
decisions, it is typically required that, for pre-specified numbers
$\de_i \in (0, 1)$,   \be \la{mainreq}
 \Pr \{ \tx{Accept} \; \mscr{H}_i \mid \se \} \geq 1 - \de_i, \qqu \fa \se
\in \varTheta_i, \qu i = 0, 1, \cd, m - 1 \ee with $\varTheta_0 = \{ \se \in \Se_0: \se \leq \se_1^{\prime} \}, \; \varTheta_{m-1} = \{ \se \in
\Se_{m-1}: \se \geq \se_{m - 1}^{\prime \prime} \}$ and $\varTheta_i = \{ \se \in \Se_i: \se_i^{\prime \prime} \leq \se \leq \se_{i+1}^{\prime}
\}, \; i = 1, \cd, m - 2$, where $\se_i^\prime, \; \se_i^{\prime \prime}$ are parametric values in ${\Se}$ such that $\se_1^\prime < \se_1, \;
\se_{m-1}^{\prime \prime} > \se_{m-1}$ and $\se_{i - 1} < \se_{i - 1}^{\prime \prime} \leq
 \se_i^\prime < \se_i < \se_i^{\prime \prime} \leq \se_{i
+ 1}^\prime < \se_{i + 1}$ for $i = 2, \cd, m - 2$. For $i = 0, 1, \cd, m - 1$, $\Pr \{ \tx{Accept} \; \mscr{H}_i \mid \se \}$ is referred to as
an Operating Characteristic (OC) function.  Since there is no requirement imposed for controlling the risk of making wrong decisions for $\se$
in ${\Se} \setminus \cup_{j = 0}^{m - 1} \varTheta_j = \cup_{i = 1}^{m-1} (\se_i^\prime, \se_i^{\prime \prime})$, such a remainder set, $\cup_{i
= 1}^{m-1} (\se_i^\prime, \se_i^{\prime \prime})$, is referred to as an {\it indifference zone}. The concept of indifference zone was introduced
by Wald \cite{Wald} for two main reasons.  First, when the parameter $\se$ is close to $\se_i$, the margin between adjacent parameter subsets
$\Se_{i-1}$ and $\Se_i$, it is immaterial to decide whether $\mscr{H}_{i-1}$ or $\mscr{H}_i$ should be accepted. Second, the sample size
required to make a reliable decision between consecutive hypotheses $\mscr{H}_{i-1}$ and $\mscr{H}_i$ becomes increasingly intolerable as $\se$
tends to $\se_i$.  Undoubtedly, the indifference zone should be sufficiently ``narrow'' so that the consequence of making erroneous decision is
practically unimportant when $\se$ lies in it. In general, a testing plan in our proposed framework consists of $s$ stages. For $\ell = 1, \cd,
s$, the number of available samples (i.e., sample size) of the $\ell$-th stage is denoted by $n_\ell$. For the $\ell$-th stage, a decision
variable $\bs{D}_\ell = \mscr{D}_\ell (X_1, \cd, X_{n_\ell})$ is defined in terms of samples $X_1, \cd, X_{n_\ell}$ such that $\bs{D}_\ell$
assumes $m+1$ possible values $0, 1, \cd, m$ with the following notion:

(i) Sampling is continued until $\bs{D}_\ell \neq 0$ for some $\ell
\in \{1, \cd, s\}$.

(ii) The hypothesis $\mscr{H}_j$ is accepted at the $\ell$-th stage
if $\bs{D}_\ell = j + 1$ and $\bs{D}_i = 0$ for $1 \leq i < \ell$.

For practical considerations, we shall only focus on sampling
schemes which are closed in the sense that $\Pr \{ \bs{D}_s = 0 \} =
0$.  For efficiency, a sampling scheme should satisfy the condition
that both  $\Pr \{ \bs{D}_1 \neq 0 \}$ and $\Pr \{ \bs{D}_{s-1} = 0
\}$ are greater than zero.

Let $\bs{l}$ denote the index of stage when the sampling is
terminated.  Then, the sample number when the sampling is
terminated, denoted by $\mbf{n}$,  is  equal to $n_{\bs{l}}$. We
shall focus on multistage sampling schemes which can be defined in
terms of estimator $\bs{\varphi}_n = \varphi(X_1, \cd, X_n)$ such
that $\bs{\varphi}_n$ is a ULE of $\se$ for every $n$ and that
$\bs{\varphi}_n$ converges in probability to $\se$ in the sense
that, for any $\vep > 0$ and $\de \in (0, 1)$, $\Pr \{ |
\bs{\varphi}_n - \se | \geq \vep \} < \de$ provided that $n$ is
sufficiently large. Such estimator $\bs{\varphi}_n$ is referred to
as a {\it Unimodal-likelihood Consistent Estimator} (ULCE) of $\se$.
Assume that, for any $\de \in (0, 1)$ and any positive number $n$, a
lower confidence limit $L( \bs{\varphi}_n, n, \de)$ and an upper
confidence limit $U( \bs{\varphi}_n, n, \de)$ can be constructed in
terms of estimator $\bs{\varphi}_n$ such that $\{ L( \bs{\varphi}_n,
n, \de) \leq \bs{\varphi}_n \leq U( \bs{\varphi}_n, n, \de) \}$ is a
sure event and that
\[
\Pr \{ \se \leq L( \bs{\varphi}_n, n, \de) \mid \se \} \leq \de,
\qqu  \Pr \{ \se \geq U( \bs{\varphi}_n, n, \de) \mid \se \} \leq
\de
\]
for any $\se \in \Se$.  To construct testing plans satisfying the requirement (\ref{mainreq}),  we propose to use fixed-sample-size confidence
limits satisfying such assumptions.   For the $\ell$-th stage, an estimator $\wh{\bs{\se}}_\ell$ for $\se$ can be defined in terms of samples
$X_1, \cd, X_{n_\ell}$ as $\bs{\varphi}_{n_\ell} = \varphi (X_1, \cd, X_{n_\ell})$.  Accordingly, the decision variables $\bs{D}_\ell$ can be
defined in terms of estimator $\wh{\bs{\se}}_\ell = \bs{\varphi}_{n_\ell}$.  The overall estimator for $\se$, denoted by $\wh{\bs{\se}}$, is
equal to $\wh{\bs{\se}}_{\bs{l}}$.

 Our general principle for constructing multistage testing plans and their properties can be
described by Theorem \ref{Multi_Comp_Exact} as follows.

 \beT
\la{Multi_Comp_Exact}

Let $\al_i = O(\ze) \in (0, 1), \; \ba_i = O(\ze) \in (0, 1)$ for $i = 1, \cd, m - 1$, where $O(\ze)$ denotes functions of the same order as
$\ze$. Let $\al_0 = \al_1, \; \ba_m = \ba_{m-1}$ and $\al_m = \ba_0 = 0$.  Let $\se_0^\prime = - \iy$ and $\se_m^{\prime \prime} = \iy$.   Let
$n_1 < n_2 < \cd < n_s$ be sample sizes such that the largest sample size $n_s$ is no less than the smallest positive integer $n$ guaranteeing
that $\{ \se_i^\prime \leq L(\bs{\varphi}_n, n, \al_i ) \leq U(\bs{\varphi}_n, n, \ba_i ) \leq \se_i^{\prime \prime} \}$ is not an impossible
event for $i = 1, \cd, m - 1$. Suppose that $\bs{\varphi}_n$ is a ULCE of $\se$ and that the decision variables $\bs{D}_\ell$ are defined such
that \bee \{ \bs{D}_\ell = i \} \subseteq \{ \se_{i-1}^\prime \leq L(\wh{\bs{\se}}_\ell, n_\ell, \al_{i-1} ) \leq U(\wh{\bs{\se}}_\ell, n_\ell,
\ba_i ) \leq \se_i^{\prime \prime} \}, \qqu 1 \leq i \leq m \eee for $\ell = 1, \cd, s$. The following statements (I)-(VI) hold true for $m \geq
2$.

(I) $\Pr \{ \tx{Reject} \; \mscr{H}_0 \mid \se \}$ is non-decreasing
with respect to $\se \in \varTheta_0$.

(II) $\Pr \{ \tx{Reject} \; \mscr{H}_{m-1} \mid \se \} \; \tx{is
non-increasing with respect to} \; \se \in \varTheta_{m-1}$.

(III)  $\Pr \{ \tx{Reject} \; \mscr{H}_i  \mid \se \} \leq s ( \max
\{ \al_j: i < j \leq m \}  + \max \{ \ba_j: 0 \leq j \leq i \}  )$
for any $\se \in \varTheta_i$ and $i = 0, 1, \cd, m - 1$.

(IV) For $0 < i \leq m- 1$, $\Pr \{ \tx{Accept} \; \mscr{H}_i \mid
\se \}$ is no greater than $s \al_i$ and is non-decreasing with
respect to $\se \in {\Se}$ no greater than $\se_i^\prime$.

(V) For $0 \leq i \leq m-2$, $\Pr \{ \tx{Accept} \; \mscr{H}_i \mid
\se \}$ is no greater than $s \ba_{i + 1}$ and is non-increasing
with respect to $\se \in {\Se}$ no less than $\se_{i+1}^{\prime
\prime}$.

(VI)  Assume that $\bb{E} [e^{\ro X}]$ exists for any $\ro \in
\bb{R}$ and that {\small $\bs{\varphi}_n = \f{ \sum_{i=1}^n X_i
}{n}$} is an unbiased and unimodal-likelihood estimator of $\se$,
where $X_1, X_2,  \cd$ are i.i.d. samples of $X$.   Then, for $i =
0, 1, \cd, m - 1$, $\lim_{\ze \to 0} \Pr \{ \tx{Reject} \;
\mscr{H}_i \mid \se \} = 0$ for any $\se \in \varTheta_i$.

 Moreover, the following statements
(VII), (VIII) and (IX) hold true for $m \geq 3$.

(VII) \bee & & \Pr \{ \tx{Reject} \; \mscr{H}_i \mid \se \} \leq \Pr
\{ \tx{Reject} \; \mscr{H}_i, \; \wh{\bs{\se}} \leq a \mid a \} +
\Pr \{ \tx{Reject} \; \mscr{H}_i, \; \wh{\bs{\se}} \geq b \mid b \},\\
&  & \Pr \{ \tx{Reject} \; \mscr{H}_i \mid \se \} \geq  \Pr \{
\tx{Reject} \; \mscr{H}_i, \; \wh{\bs{\se}} \leq a \mid b \} + \Pr
\{ \tx{Reject} \; \mscr{H}_i, \; \wh{\bs{\se}} \geq b \mid a \} \eee
for any $\se \in [a, b ] \subseteq \varTheta_i$ and $1 \leq i \leq m
- 2$.

(VIII) $\Pr \{ \tx{Reject} \; \mscr{H}_0 \; \tx{and} \;
\mscr{H}_{m-1} \mid \se \}$ is non-decreasing with respect to $\se
\in \varTheta_0$ and is non-increasing with respect to $\se \in
\varTheta_{m-1}$.

(IX) $\Pr \{ \tx{Reject} \; \mscr{H}_0 \; \tx{and} \; \mscr{H}_{m-1}
\mid \se \}$ is no greater than $s \times \max \{ \al_i: 1 \leq i
\leq m - 2 \}$ for $\se \in \varTheta_0$ and is no greater than $s
\times \max \{ \ba_i: 2 \leq i \leq m - 1 \}$ for $\se \in
\varTheta_{m-1}$.

\eeT

Theorem \ref{Multi_Comp_Exact} asserts that the probabilities of committing decision errors for the proposed testing plan can be adjusted below
any prescribed level by choosing a sufficiently small value of $\ze > 0$.  This process of adjusting the probabilities of committing decision
errors can be accomplished by the bisection risk tuning techniques established in our earlier paper \cite{Chen_TH}.  Actually, we had derived in
\cite{Chen_TH} many concrete testing plans by the same methodology as that of the general stopping and decision rules described by Theorem
\ref{Multi_Comp_Exact}.  To see the fundamental principle behind the construction of stopping and decision rules, we will demonstrate that the
general problem of designing multistage plans for testing multiple hypotheses can be reformulated as the problem of constructing a {\it
sequential random interval} to satisfy certain specifications of coverage probability. Note that the word ``sequential'' is used to indicate the
fact that the random interval is constructed from samples of random size. As a consequence of the connection between testing plans and
sequential random intervals, the properties of the testing plan described by statements (I)--(IX) of Theorem \ref{Multi_Comp_Exact} can be
justified by the general theory of coverage probability of sequential random intervals established in Sections 2.5 and 2.6 of \cite{Chen_EST}.

To illustrate why the general hypothesis testing problem defined by (\ref{mainpr}) and (\ref{mainreq}) can be cast into the framework of
constructing sequential random intervals with pre-specified coverage probabilities, consider a multistage sampling scheme with sample sizes $n_1
< n_2 < \cd < n_s$ for random variable $X$ parameterized by $\se \in \Se$. Let $\wh{\bs{\se}} = \varphi(X_1, \cd, X_{\mbf{n}})$, where $\mbf{n}$
is the number of samples at the termination of the process of drawing samples $X_1, X_2, \cd$. Let $\mscr{L} ( \wh{\bs{\se}}, \mbf{n} )$ and
$\mscr{U} ( \wh{\bs{\se}}, \mbf{n} )$ be functions of $\wh{\bs{\se}}$ and $\mbf{n}$ such that the following requirements are satisfied:

(i): The sequential random interval {\small $( \mscr{L} ( \wh{\bs{\se}}, \mbf{n} ), \mscr{U} ( \wh{\bs{\se}}, \mbf{n} ) )$} has $m$ possible
realizations: {\small $( \se_i^\prime, \se_{i+1}^{\prime \prime} ), \;  i = 0,  1, \cd, m - 1$}.

(ii):  \be \la{coveragereq} \Pr \{ \mscr{L} ( \wh{\bs{\se}}, \mbf{n} ) < \se < \mscr{U} ( \wh{\bs{\se}}, \mbf{n} ) \mid \se \} > 1 - \de_i \ee
for any $\se \in \varTheta_i$ and $i = 0, 1, \cd, m - 1$.

\bsk

Given that the sequential random interval $( \mscr{L} ( \wh{\bs{\se}}, \mbf{n} ), \; \mscr{U} ( \wh{\bs{\se}}, \mbf{n} ))$ satisfying
requirements (i)--(ii) is constructed, the risk requirement (\ref{mainreq}) can be satisfied by using the sequential random interval $( \mscr{L}
( \wh{\bs{\se}}, \mbf{n} ), \; \mscr{U} ( \wh{\bs{\se}}, \mbf{n} ))$ to define a decision rule such that, for $i = 0, 1, \cd, m - 1$, hypothesis
$\mscr{H}_i$ is accepted when the sequential random interval $( \mscr{L} ( \wh{\bs{\se}}, \mbf{n} ), \mscr{U} ( \wh{\bs{\se}}, \mbf{n} ) )$
takes the $i$-th realization at the termination of the sampling process.   Our claim that the risk requirement (\ref{mainreq}) is satisfied can
be justified by virtue of (\ref{coveragereq}) and the observation that \be \la{have8} \{ \tx{Accept} \; \mscr{H}_i \} =\{ \mscr{L} (
\wh{\bs{\se}}, \mbf{n} ) < \se < \mscr{U} ( \wh{\bs{\se}}, \mbf{n} ) \} \ee for any $\se \in \varTheta_i$ and $i = 0, 1, \cd, m - 1$.

The construction of the sequential random interval {\small $( \mscr{L} ( \wh{\bs{\se}}, \mbf{n} ), \mscr{U} ( \wh{\bs{\se}}, \mbf{n} ) )$}
fulfilling requirements (i)--(ii) can be accomplished as follows. Noting that
\[
\Pr \{ \mscr{L} ( \wh{\bs{\se}}, \mbf{n} ) < \se < \mscr{U} (
\wh{\bs{\se}}, \mbf{n} ) \mid \se \} = \Pr \{ ( \mscr{L} (
\wh{\bs{\se}}, \mbf{n} ), \mscr{U} ( \wh{\bs{\se}}, \mbf{n} ) )
\subseteq (\se_i^\prime, \se_{i+1}^{\prime \prime} ) \mid \se \}
\]
for any $\se \in \varTheta_i$ with $i \in \{0, 1, \cd, m - 1\}$, we propose to control the coverage probability of the sequential random
interval by creating a particular inclusion relationship between $m$ sequences of confidence intervals and the desired sequential random
interval. For $i = 0, 1, \cd, m - 1$, let $( \se_i^\prime, \se_{i+1}^{\prime \prime} )$ be referred to as the $i$-th realization of the
sequential random interval $( \mscr{L} ( \wh{\bs{\se}}, \mbf{n} ), \mscr{U} ( \wh{\bs{\se}}, \mbf{n} ) )$.  For $\ell = 1, \cd, s$, let
$\wh{\bs{\se}}_\ell = \bs{\varphi}_{n_\ell} = \varphi (X_1, \cd, X_{n_\ell})$ as before.  For $i = 0, 1, \cd, m - 1$, let
$(L(\wh{\bs{\se}}_\ell, n_\ell, \al_i ), U(\wh{\bs{\se}}_\ell, n_\ell, \ba_{i+1} )), \; \ell = 1, \cd, s$ be referred to as the $i$-th {\it
controlling confidence sequence}.   Our purpose of introducing the term ``controlling confidence sequence'' are two folded: (i) To indicate that
the confidence sequences are used to control the coverage probability of the desired sequential random interval; (ii) To avoid the confusion
between the sequences of confidence intervals and the sequential random intervals.  In this setting, we propose the following stopping and
decision rule:

{\it Continue the sampling process if there exists no index $i\in \{0, 1, \cd, m - 1\}$ such that the $i$-th controlling confidence sequence is
included by the $i$-th realization of the sequential random interval.  At the termination of the sampling process, make the following decision:

(a):  If there exists a unique index $i \in \{0, 1, \cd, m - 1\}$ such that the $i$-th controlling confidence sequence is included by the $i$-th
realization of the sequential random interval, then designate the $i$-th realization as the outcome of the sequential random interval.

(b):  If there exist two consecutive indexes $i - 1$ and $i$ in $\{1, \cd, m - 1\}$ such that the $(i - 1)$-th and $i$-th controlling confidence
sequences are included, respectively, by the $(i-1)$-th and $i$-th realizations of the sequential random interval, then designate either the
$(i-1)$-th or $i$-th realization as the outcome of the sequential random interval based on an arbitrarily pre-specified policy. } \bsk

As a consequence of the assumption that all controlling confidence sequences always include $\wh{\bs{\se}}_\ell$ for $\ell = 1, \cd, s$, it is
impossible that there are more than two indexes $i\in \{0, 1, \cd, m - 1\}$ such that the $i$-th controlling confidence sequence is included by
the $i$-th realization of the sequential random interval.   However, the situation described in (b) is possible.   One possible policy to handle
such a tied situation for the $\ell$-th stage is to designate the $(i-1)$-th realization as the outcome of the sequential random interval if
 $\wh{\bs{\se}}_\ell \leq m_{\ell, i}$ and designate the $i$-th
 realization as the outcome of the sequential random interval if
 $\wh{\bs{\se}}_\ell > m_{\ell, i}$, where
 $m_{\ell, i}$ can be taken as $m_{\ell, i} = \f{ \se_i^{\prime}  +  \se_i^{\prime \prime} }{2}$ or $m_{\ell, i}  =
 \f{1}{2} [  \min \{ z \in I_{\wh{\bs{\se}}_\ell}: L(z, n_\ell,
\al_i ) \geq \se_i^{\prime} \} +  \max \{ z \in I_{\wh{\bs{\se}}_\ell}: U(z, n_\ell, \ba_i ) \leq \se_i^{\prime \prime} \} ]$ for $\ell = 1,
\cd, s$ and $i = 1, \cd, m - 1$.

Based on the above stopping and decision rule, the decision variables can be expressed as in terms of $\wh{\bs{\se}}_\ell$.  For this purpose,
define \bee & & \mscr{A}_{\ell, i} = \{ z \in I_{\wh{\bs{\se}}_\ell}: L (z, n_\ell, \al_i )
\geq \se_i^{\prime} \},  \\
&  & \mscr{B}_{\ell, i} = \{ z \in I_{\wh{\bs{\se}}_\ell}: U (z, n_\ell, \ba_i ) \leq \se_i^{\prime \prime} \},\\
&  & \mscr{C}_{\ell, i} = \{ z \in I_{\wh{\bs{\se}}_\ell}: \se_i^{\prime} \leq L (z, n_\ell, \al_i ) \leq U
(z, n_\ell, \ba_i ) \leq \se_i^{\prime \prime} \},\\
&  & f_{\ell, i} = \bec \max \mscr{B}_{\ell, i} & \tx{if} \;
\mscr{B}_{\ell, i} \neq \emptyset \; \tx{and} \; \mscr{C}_{\ell, i} = \emptyset,\\
c_{\ell, i}  & \tx{if} \; \mscr{C}_{\ell, i} \neq
\emptyset,\\
- \iy & \tx{if} \; \mscr{B}_{\ell, i} = \emptyset \eec\\
&  & g_{\ell, i} = \bec \min \mscr{A}_{\ell, i} & \tx{if} \;
\mscr{A}_{\ell, i} \neq \emptyset \; \tx{and} \; \mscr{C}_{\ell, i} = \emptyset,\\
c_{\ell, i}  & \tx{if} \; \mscr{C}_{\ell, i} \neq
\emptyset,\\
\iy & \tx{if} \; \mscr{A}_{\ell, i} = \emptyset \eec \eee for $\ell = 1, \cd, s$ and $i = 1, \cd, m - 1$, where $c_{\ell, i}$ is a number such
that $\min \mscr{C}_{\ell, i} \leq c_{\ell, i} \leq \max \mscr{C}_{\ell, i}$.  In particular, we can take $c_{\ell, i} = \f{1}{2} (  \min
\mscr{C}_{\ell, i} + \max \mscr{C}_{\ell, i}  )$ or $c_{\ell, i} = \f{1}{2} ( \se_i^\prime + \se_i^{\prime \prime} )$ for $i = 1, \cd, m - 1$.
Define $f_{\ell, m} = \iy$ and $g_{\ell, 0} = - \iy$.  Then, the decision variables can be defined as  {\small \be \la{defgoog} \bs{D}_\ell =
\bec i & \tx{if there exists} \; i \in \{1, \cd, m \} \; \tx{such that}
\; \; g_{\ell, i-1} < \wh{\bs{\se}}_\ell \leq f_{\ell,i}\\
0 & \tx{otherwise} \eec \ee} for $\ell = 1, \cd, s$.  At the termination of the sampling process, the $(i-1)$-th
 realization is designated as the outcome of the sequential random interval $( \mscr{L} ( \wh{\bs{\se}}, \mbf{n} ), \;
\mscr{U} ( \wh{\bs{\se}}, \mbf{n} ))$ if $\bs{D}_\ell = i$.

Under the assumption that $\al_i = O(\ze)$ and $\ba_i = O(\ze)$ for $i = 1, \cd, m - 1$ and that the maximum sample size $n_s$ is no less than
the minimum integer $n$ guaranteeing that $\{ \se_i^\prime \leq L(\bs{\varphi}_n, n, \al_i ) \leq U(\bs{\varphi}_n, n, \ba_i ) \leq
\se_i^{\prime \prime} \}$ is not an impossible event for $i = 1, \cd, m - 1$, we can show that the coverage requirement (\ref{coveragereq}) can
be satisfied with the above stopping and decision rule by choosing $\ze > 0$ to be a sufficiently small number.

From above discussion on the reformulation of testing procedures, it can be seen that the testing plan described in Theorem
\ref{Multi_Comp_Exact} is proposed by the same idea of constructing the sequential random interval as above.  Accordingly, the properties of the
testing plan can be shown by virtue of (\ref{have8}) and the theory of coverage probability of sequential random intervals established in
Sections 2.5 and 2.6 of \cite{Chen_EST}.  Most importantly,  such reformulation demonstrates that the two seemingly different fields of
sequential analysis, multistage parameter estimation and hypothesis testing, can be cast into a much broader and unified framework of
constructing sequential random intervals of pre-specified coverage probabilities.  A critical methodology to make such unification possible is
the use of multiple confidence sequences to control the coverage probability of the desired sequential random interval. The controllability of
the coverage probability is due to the inclusion relationship implemented by the stopping and decision rules.  In view of its versatility, we
call such a methodology of using confidence sequences to define stopping and decision rules as the Inclusion Principle (IP). In addition to IP,
the design of multistage testing plans and estimation procedures can be tackled by the same computational method -- bisection coverage tuning,
which is discussed in our papers \cite{Chen_EST} and \cite{Chen_TH}.

Haven established the fundamental principle of constructing testing plans, our next task is to make the stopping and decision rules as simple as
possible. In situations that the parameter $\se$ to be tested is the expectation of $X$, for the sake of simplicity, one can replace the strict
confidence limits by their approximations derived from the central limit theorem or by their bounds derived from probabilistic inequalities such
as Chernoff bound.  In particular, we can use the following methods for constructing approximate confidence intervals.  Assume that $X_1, X_2,
\cd$ are identical samples of $X$ and that the variance of the sample mean $\ovl{X}_n \DEF \f{\sum_{i=1}^n X_i}{n}$ is a bivariate function,
denoted by $\mscr{V} (\se, n)$, of $\se$ and $n$. If the sample size $n$ is large, then the central limit theorem  may be applied to establish
the normal approximation \bel & & F_{\ovl{X}_n} \li ( z, \se \ri ) \DEF \Pr \{ \ovl{X}_n \leq z \mid \se \} \ap \Phi \li ( \f{ z - \se  } {
\sqrt{ \mscr{V} ( \se, n )
} } \ri ), \la{normap1}\\
&  & G_{\ovl{X}_n} \li ( z, \se  \ri ) \DEF \Pr \{ \ovl{X}_n \geq z \mid \se  \} \ap \Phi \li ( \f{ \se - z } { \sqrt{ \mscr{V} ( \se, n) } }
\ri ). \la{normap2} \eel Let $\de \in (0, 1)$. If $\ovl{X}_n$ is a ULE of $\se$, then we can obtain lower and upper confidence limits
respectively as
\[
L (\ovl{X}_n, n, \de) = \inf \li \{  \se \in \Se : \Phi \li ( \f{ \se - \ovl{X}_n } { \sqrt{ \mscr{V} ( \se, n) } } \ri ) > \f{\de}{2} \ri \}
\]
and
\[
U (\ovl{X}_n, n, \de) = \sup \li \{  \se \in \Se : \Phi \li ( \f{ \ovl{X}_n - \se  } { \sqrt{ \mscr{V} ( \se, n ) } } \ri ) > \f{\de}{2} \ri \}
\]
such that \[ \Pr \{ L (\ovl{X}_n, n, \de) \geq \se \mid \se \} \ap \f{\de}{2}, \qqu  \Pr \{ U (\ovl{X}_n, n, \de) \leq \se \mid \se \} \ap
\f{\de}{2}
\]
and $\Pr \{ L (\ovl{X}_n, n, \de) < \se < U (\ovl{X}_n, n, \de) \mid \se \} \ap 1 - \de$ for $\se \in \Se$.  To improve the accuracy of normal
approximation (\ref{normap1}) and (\ref{normap2}), we propose to replace $\se$ in $\mscr{V} ( \se, n )$ by $z + w (\se - z)$ with $w \in [0,
1]$. That is, we suggest modifying (\ref{normap1}) and (\ref{normap2}) as follows: \bel & & F_{\ovl{X}_n} \li ( z, \se \ri ) \ap
\Phi \li ( \f{ z - \se  } { \sqrt{ \mscr{V} ( z + w (\se - z), n ) } } \ri ), \la{normap3}\\
&  & G_{\ovl{X}_n} \li ( z, \se  \ri ) \ap \Phi \li ( \f{ \se - z } { \sqrt{ \mscr{V} ( z + w (\se - z), n) } } \ri ). \la{normap4} \eel The new
parameter $w$ is introduced to improve the accuracy of approximation.   Based on the new approximation (\ref{normap3}) and (\ref{normap4}), we
propose to obtain lower and upper confidence limits for $\se$ respectively as \be \la{best_CI_A}
 L (\ovl{X}_n, n, \de) = \inf
\li \{ \se \in \Se : \Phi \li ( \f{ \se - \ovl{X}_n } { \sqrt{ \mscr{V} ( \ovl{X}_n + w (\se - \ovl{X}_n), n) } } \ri ) > \f{\de}{2} \ri \} \ee
and \be \la{best_CI_B} U (\ovl{X}_n, n, \de) = \sup \li \{  \se \in \Se : \Phi \li ( \f{ \ovl{X}_n - \se  } { \sqrt{ \mscr{V} ( \ovl{X}_n + w
(\se - \ovl{X}_n), n ) } } \ri ) > \f{\de}{2} \ri \} \ee so that $\Pr \{ L (\ovl{X}_n, n, \de) \geq \se \mid \se \} \ap \f{\de}{2}, \; \Pr \{ U
(\ovl{X}_n, n, \de) \leq \se \mid \se \} \ap \f{\de}{2}$ and $\Pr \{ L (\ovl{X}_n, n, \de) < \se < U (\ovl{X}_n, n, \de) \mid \se \} \ap 1 -
\de$ for $\se \in \Se$.  The confidence limits constructed by this approach can be used to derive simple stopping and decision rules. Moreover,
the performance of the testing plans can be optimized with respect to $w \in [0, 1]$. It should be noted that for parameters of binomial,
Poisson, geometric and hypergeometric distributions, explicit forms of confidence limits can be obtained from (\ref{best_CI_A}) and
(\ref{best_CI_B}) by solving quadratic equations. For example, in scenarios that  $X_1, \cd, X_n$ are i.i.d. samples of Bernoulli random
variable $X$ such that $\Pr \{ X = 1 \} = 1 - \Pr \{ X = 0 \} = p \in (0, 1)$, the lower and upper confidence limits for $p$ can be readily
derived from (\ref{best_CI_A}) and (\ref{best_CI_B}) respectively by
\[
L (\ovl{X}_n, n, \de) = \f{ \ovl{X}_n +  \f{ w \mcal{Z}^2 }{ 2 n} [ 1 - 2 (1 - w) \ovl{X}_n ]   - \mcal{Z} \sq{ \f{ \ovl{X}_n ( 1 - \ovl{X}_n )
}{n} + \li ( \f{ w \mcal{Z}  }{ 2 n } \ri )^2 } }{ 1 + \f{ ( w \mcal{Z} )^2 }{n}  }
\]
and
\[
U (\ovl{X}_n, n, \de) =  \f{ \ovl{X}_n  +  \f{ w \mcal{Z}^2 }{ 2 n} [ 1 - 2 (1 - w) \ovl{X}_n ] + \mcal{Z} \sq{ \f{ \ovl{X}_n ( 1 - \ovl{X}_n )
}{n} + \li ( \f{ w \mcal{Z}  }{ 2 n } \ri )^2 } }{ 1 + \f{ ( w \mcal{Z} )^2 }{n}  },
\]
where $\mcal{Z}$ is the critical value such that $\Phi( \mcal{Z} ) = 1 - \f{\de}{2}$.

\sect{Multistage Tests for Functions of Two Binomial Proportions}

Testing hypotheses on functions of two binomial proportions is particularly important in prospective comparative studies such as randomized
controlled clinical trial.  More formally, let $X$ and $Y$ be independent Bernoulli random variables such that $\Pr \{ X = 1 \} = 1 - \Pr \{ X =
0 \} = p_x \in (0, 1)$ and $\Pr \{ Y = 1 \} = 1 - \Pr \{ Y = 0 \} = p_y \in (0, 1)$.  Let $g(.,.)$ be a bivariate function of $p_x$ and $p_y$.
It is a frequent problem to test hypotheses regarding $g(p_x, p_y)$ based on samples of $X$ and $Y$. Typical examples of $g(p_x, p_y)$ are
$g(p_x, p_y) = p_x - p_y$ and $g(p_x, p_y) = \f{p_x}{p_y}$, which are respectively referred to as the {\it difference of population proportions}
and {\it ratio of population proportions}.  In the sequel, we shall design multistage plans for testing hypotheses regarding $g(p_x, p_y)$ by
virtue of the method proposed in Section 5.  The main idea is to make use of confidence limits to define stopping and decision rules. To address
the problems of testing hypotheses regarding $g(p_x, p_y)$ in a unified framework, we define $\se = g(p_x, p_y)$ and consider the general
problem of testing $m$ mutually exclusive and exhaustive composite hypotheses defined by (\ref{mainpr}) and (\ref{mainreq}).  For simplicity of
notations, let $\bs{p} = (p_x, p_y)$ and abbreviate $g(p_x, p_y)$ as $g(\bs{p})$.

Our strategy is to seek a class of testing plans such that the risk of making wrong decisions can be controlled by the risk tuning parameter
$\ze > 0$, and then apply the risk tuning technique to adjust the risk level to satisfy the requirement (\ref{mainreq}).  Let $O(\ze)$ denote
the class of functions of the same order of $\ze$.  Let $\al_i = O(\ze) \in (0, 1)$ and $\ba_i = O(\ze) \in (0, 1)$ for $i = 0, 1, \cd, m$.  Let
$\se_0^\prime = - \iy$ and $\se_m^{\prime \prime} = \iy$ as before. Specifically, we seek a class of sampling schemes associated with $\ze$
satisfying the following requirements:

\bed

\item  (i) The sample sizes $N_{1, x} (\ze) < N_{2, x} (\ze)  < \cd < N_{s, x} (\ze)$ and $N_{1, y} (\ze)  < N_{2, y} (\ze) < \cd < N_{s, y}
(\ze)$ are decreasing functions of the risk tuning parameter $\ze > 0$. Namely, the sample sizes are increasing as $\ze$ decreases. For
simplicity of notations, we abbreviate $N_{\ell, x}(\ze), \; N_{\ell, y}(\ze)$ as $N_{\ell, x}, \; N_{\ell, y}$ respectively and  define \bee &
&  K_{\ell, x} = \sum_{i = 1}^{N_{\ell, x}} X_i, \qqu K_{\ell, y} = \sum_{i = 1}^{N_{\ell, y}} Y_i, \qqu
\bs{K}_\ell = (K_{\ell, x}, K_{\ell, y}), \\
&  & \wh{\bs{p}}_{\ell, x} = \f{K_{\ell, x}}{N_{\ell, x}}, \qqu \qu \wh{\bs{p}}_{\ell, y} = \f{K_{\ell, y}}{N_{\ell, y}}, \qqu \qu
\wh{\bs{p}}_{\ell} = (\wh{\bs{p}}_{\ell, x}, \wh{\bs{p}}_{\ell, y}) \eee for $\ell = 1, \cd, s$.

\item (ii) For $\ell = 1, \cd, s$, one-sided  confidence intervals $[\mcal{L}(\wh{\bs{p}}_{\ell}, \ell, \al_i), \; \iy )$ and $(- \iy, \;
\mcal{U}(\wh{\bs{p}}_{\ell}, \ell, \ba_i)]$ for $\se = g(\bs{p})$ can be constructed in terms of estimators $\wh{\bs{p}}_{\ell}$, sample sizes
$N_{\ell, x}, \; N_{\ell, y}$ and the risk tuning parameter $\ze
> 0$ such that \bee &  & \mcal{L}(\wh{\bs{p}}_{\ell}, \ell, \al_i) \leq g( \wh{\bs{p}}_{\ell} ),
\qqu \Pr \{ \mcal{L}(\wh{\bs{p}}_{\ell}, \ell, \al_i) \geq
\se \mid \bs{p} \} \leq \al_i,\\
&  & \mcal{U}(\wh{\bs{p}}_{\ell}, \ell, \ba_i) \geq g( \wh{\bs{p}}_{\ell} ), \qqu \Pr \{ \mcal{U}(\wh{\bs{p}}_{\ell}, \ell, \ba_i) \leq \se \mid
\bs{p} \} \leq \ba_i, \eee for $i = 0, 1, \cd,  m$ and $\bs{p} \in \{ (p_x, p_y): p_x \in (0, 1), \; p_y \in (0, 1) \}$.

\item (iii) The sampling process is continued until $\se_{i-1}^\prime \leq \mcal{L} (\wh{\bs{p}}_\ell, \ell, \al_{i-1} ) \leq \mcal{U}
(\wh{\bs{p}}_\ell, \ell, \ba_i ) \leq \se_i^{\prime \prime}$ for some $i \in \{1, \cd, m \}$ at some stage with index $\ell \in \{1, \cd, s \}$.

\item (iv) When the sampling process is terminated at the $\ell$-th stage, there are only two possible cases (A) and (B) described as follows:

  \bed
  \item [Case (A)]: There is a unique index $i \in \{1, \cd, m\}$ such that $\{ \se_{i-1}^\prime \leq \mcal{L} (\wh{\bs{p}}_\ell, \ell, \al_{i-1} ) \leq \mcal{U}
(\wh{\bs{p}}_\ell, \ell, \ba_i ) \leq \se_i^{\prime \prime} \}$.

 \item [Case (B)]: There is a unique index $i \in \{1, \cd, m - 1\}$ such that $\{ \se_{i-1}^\prime \leq \mcal{L} (\wh{\bs{p}}_\ell, \ell, \al_{i-1} ) \leq \mcal{U}
(\wh{\bs{p}}_\ell, \ell, \ba_i ) \leq \se_i^{\prime \prime} \}$ and $\{ \se_{i}^\prime \leq \mcal{L} (\wh{\bs{p}}_\ell, \ell, \al_{i} ) \leq
\mcal{U} (\wh{\bs{p}}_\ell, \ell, \ba_{i+1} ) \leq \se_{i+1}^{\prime \prime} \}$.

\item The decision is made based on the following policy:

 In Case (A), accept hypothesis $\mscr{H}_{i-1}$.  In Case (B), accept $\mscr{H}_{i-1}$ if $g( \wh{\bs{p}}_\ell ) \leq \f{
\se_{i}^\prime +  \se_i^{\prime \prime} }{2}$ and accept $\mscr{H}_{i}$ otherwise.

\eed

\item (v) The sampling process is guaranteed to be terminated at or before the $s$-th stage.

\eed

\bsk

Within this class of sampling schemes, we can apply the bisection risk tuning technique to determine $\ze > 0$ as large as possible such that
the risk requirement (\ref{mainreq}) is satisfied.   In practices, the inequalities ``$\Pr \{ \mcal{L}(\wh{\bs{p}}_{\ell}, \ell, \al_i) \geq \se
\mid \bs{p} \} \leq \al_i$'' and ``$\Pr \{ \mcal{U}(\wh{\bs{p}}_{\ell}, \ell, \ba_i) \leq \se \mid \bs{p} \} \leq \ba_i$'' in requirement (ii)
can be respectively relaxed as
\[
\Pr \{ \mcal{L}(\wh{\bs{p}}_{\ell}, \ell, \al_i) \geq \se \mid \bs{p} \} \ap \al_i,  \qqu  \Pr \{ \mcal{U}(\wh{\bs{p}}_{\ell}, \ell, \ba_i) \leq
\se \mid \bs{p} \} \ap \ba_i
\]
provided that the risk of making wrong decisions can be controlled by $\ze$.   The constraints ``$\mcal{L}(\wh{\bs{p}}_{\ell}, \ell, \al_i) \leq
g( \wh{\bs{p}}_{\ell} )$'' and ``$\mcal{U}(\wh{\bs{p}}_{\ell}, \ell, \ba_i) \geq g( \wh{\bs{p}}_{\ell} )$'' in requirement (ii) ensure that
``Case (A)'' and ``Case (B)'' are the only two possible cases at the termination of the sampling process.

As an illustration of the construction of sampling schemes, we consider the problem defined by (\ref{mainpr}) and (\ref{mainreq}) of testing $m$
mutually exclusive and exhaustive composite hypotheses regarding $\se = g (\bs{p}) = p_x - p_y$.   At each stage, one can use the following
confidence limits {\small \bee \mcal{L}(\wh{\bs{p}}_{\ell}, \ell, \al_i) = \mscr{L} ( \wh{\bs{p}}_{\ell, x}, \wh{\bs{p}}_{\ell, y}, N_{\ell, x},
N_{\ell, y}, \al_i ), \qqu \mcal{U}(\wh{\bs{p}}_{\ell}, \ell, \ba_i) = \mscr{U} ( \wh{\bs{p}}_{\ell, x}, \wh{\bs{p}}_{\ell, y}, N_{\ell, x},
N_{\ell, y}, \ba_i ), \qu i = 0, 1, \cd, m \eee} for $\ell = 1, \cd, s$, where $\mscr{L}$ and $\mscr{U}$ are confidence limits investigated in
\cite{Newcombe}, which can be expressed as \bee &  & \mscr{L} ( \wh{p}_x, \wh{p}_y, N_x, N_y, \de ) = \wh{p}_x - \wh{p}_y - \mcal{Z}_{\de} \sq{
\f{ l_x (1 - l_x) }{N_x} +  \f{ u_y (1 - u_y) }{N_y} }, \\
&  & \mscr{U} ( \wh{p}_x, \wh{p}_y, N_x, N_y, \de ) = \wh{p}_x - \wh{p}_y + \mcal{Z}_{ \de } \sq{  \f{ u_x (1 - u_x) }{N_x} +  \f{ l_y (1 - l_y)
}{N_y} }, \eee where \bee &  & \wh{p}_x \in \li \{ \f{k}{N_x}: k = 0, 1, \cd, N_x  \ri \}, \qqu \qqu \qqu \qu \; \;
\wh{p}_y \in \li \{ \f{k}{N_y}: k = 0, 1, \cd, N_y  \ri \},\\
 &  &  l_x = \f{ c + 2 N_x \wh{p}_x - \sq{ c^2 + 4 c N_x \wh{p}_x (1 - \wh{p}_x)    }  } { 2 (c + N_x) },
\qqu \qu u_x = \f{ c + 2 N_x \wh{p}_x + \sq{ c^2 + 4 c N_x \wh{p}_x (1 - \wh{p}_x) }  } { 2 (c + N_x) },\\
&  & l_y = \f{ c + 2 N_y \wh{p}_y - \sq{ c^2 + 4 c N_y \wh{p}_y (1 - \wh{p}_y)    }  } { 2 (c + N_y) }, \qqu \qu u_y = \f{ c + 2 N_y \wh{p}_y +
\sq{ c^2 + 4 c N_y \wh{p}_y (1 - \wh{p}_y) }  } { 2 (c + N_y) } \eee with $c = \mcal{Z}_{\de}^2$. It can be checked that $l_x, \; u_x$  are the
roots for $p$ in the quadratic equation $| \wh{p}_x - p | = \mcal{Z}_{ \de } \sq{ p (1 - p) \sh N_x }$ and that $l_y, \; u_y$ are the roots for
$p$ in the quadratic equation $| \wh{p}_y - p | = \mcal{Z}_{\de} \sq{ p (1 - p) \sh N_y}$.  Actually, there are many methods to construct
confidence intervals satisfying requirement (i) (see, e.g., \cite{Brown, Newcombe} and the references therein).

For a given value of the risk tuning parameter $\ze$, we have to determine sample sizes $(N_{\ell, x}, N_{\ell, y})$, $\ell = 1, \cd, s$ such
that the sampling process is guaranteed to be terminated at or before the $s$-th stage.  This can be done as follows.

In applications, the sample sizes for $Y$ can be chosen to be increasing functions of sample sizes of $X$.  In other words, one can define an
increasing function $\digamma (.)$ such that $N_{\ell, y} = \digamma (N_{\ell, x})$ for $\ell = 1, \cd, s$.  Therefore, for a given value of the
risk tuning parameter $\ze$, the maximum sample sizes $(N_{s, x}, N_{s, y})$ can be determined as $(N_x, N_y)$, where $N_y = \digamma (N_x)$
with $N_x$ guaranteeing that $\{ (\wh{p}_x, \wh{p}_y): \se_{i}^\prime \leq \mscr{L} ( \wh{p}_x, \wh{p}_y, N_x, N_y, \al_{i} ) \leq \mscr{U} (
\wh{p}_x, \wh{p}_y, N_x, N_y, \ba_i ) \leq \se_i^{\prime \prime} \}$ is nonempty for $i = 1, \cd, m - 1$.  The minimum sample sizes $(N_{1, x},
N_{1, y})$ can be determined as $(N_x, N_y)$, where $N_y = \digamma (N_x)$ with $N_x$  guaranteeing that
\[
\bigcup_{i = 1}^{m} \{ (\wh{p}_x, \wh{p}_y): \se_{i-1}^\prime \leq \mscr{L} ( \wh{p}_x, \wh{p}_y, N_x, N_y, \al_{i-1} ) \leq \mscr{U} (
\wh{p}_x, \wh{p}_y, N_x, N_y, \ba_i ) \leq \se_i^{\prime \prime} \} \] is nonempty.  To reduce the sampling cost, one should choose $N_x$ as
small as possible in the determination of $(N_{s, x}, N_{s, y})$ and $(N_{1, x}, N_{1, y})$.  Once the sample sizes for the first and last
stages are determined, one can choose the sample sizes for other stages.  For example, the sample sizes can be chosen as arithmetic or geometric
progressions.

Given that a class of testing plans associated with $\ze$ have been constructed such that the risk of making wrong decisions can be made
arbitrarily small if $\ze$ is sufficiently small, the bisection risk tuning technique can be applied to determine $\ze$ as large as possible
such that the risk requirement (\ref{mainreq}) can be satisfied.  A critical subroutine of the bisection risk tuning method is to apply Adapted
Branch and Bound algorithm proposed in Appendix T of \cite{Chen_EST} to determine whether the probability of rejecting $\mscr{H}_i$ is no
greater  than $\de_i$ for all $\bs{p}$ such that $\se = g(\bs{p}) \in \varTheta_i$. This requires bounding the probability $\Pr \{
\tx{Rejecting} \; \mscr{H}_i \mid \bs{p} \}$ for $\bs{p}$ in a rectangular domain $\mcal{Q} \DEF \{ (p_x, p_y) : 0 \leq \udl{p}_x \leq p_x \leq
\ovl{p}_x \leq 1, \; 0 \leq \udl{p}_y \leq p_y \leq \ovl{p}_y \leq 1 \}$.  From our proposed method of designing multistage testing plans, it
can be seen that, for a given value of risking tuning parameter $\ze > 0$, a testing plan can be described as follows:

At the $\ell$-th stage, the sampling process is continued if $\wh{\bs{p}}_\ell \in \mscr{C}_\ell$. Otherwise, the sampling process is terminated
and hypothesis $\mscr{H}_i$ is accepted if $\wh{\bs{p}}_\ell \in \mscr{D}_{i, \ell}$.  Here $\mscr{C}_\ell$ and $\mscr{D}_{i, \ell}, \; i = 0,
1, \cd, m-1$ are mutually exclusive and exhaustive subsets of the support of $\wh{\bs{p}}_\ell$.   In other words, the sampling plan is to
continue sampling until $\wh{\bs{p}}_\ell \in \mscr{D}_{i, \ell}$ for some $\ell \in \{1, \cd, s\}$ and $i \in \{0, 1, \cd, m - 1 \}$ and then
accept hypothesis $\mscr{H}_i$.

With the above understanding of the testing plan, one can see that for $i = 0, 1, \cd, m - 1$, \be \la{vip898} \{ \tx{Rejecting} \; \mscr{H}_i
\}  = \cup_{j = 1}^s \li \{ \wh{\bs{p}}_\ell \in \mscr{C}_\ell, \; 1 \leq \ell < j; \;  \wh{\bs{p}}_i \in \mscr{D}_{i, j}^r \ri \}, \ee where
\[
\mscr{D}_{i, j}^r \DEF \bigcup_{ 0 \leq k < m \atop{k \neq i}} \mscr{D}_{k, j}
\]
is the rejection region of hypothesis $\mscr{H}_i$.  To reduce the computational complexity for bounding $\Pr \{ \tx{Rejecting} \; \mscr{H}_i
\mid \bs{p} \}$, we invoke the truncation method proposed in \cite{ChenTr}.   For this purpose, let $\eta \in (0, 1)$ and define {\small \bee &
& a_{\ell, x} = \mrm{T}_{\mrm{lb}} (\udl{p}_x, N_{\ell, x}, \eta),  \qu b_{\ell, x} = \mrm{T}_{\mrm{ub}} (\ovl{p}_x, N_{\ell, x}, \eta),
\qu a_{\ell, y}  = \mrm{T}_{\mrm{lb}} (\udl{p}_y, N_{\ell, y}, \eta),  \qu  b_{\ell, y} = \mrm{T}_{\mrm{ub}} (\ovl{p}_y, N_{\ell, y}, \eta), \\
&  & c_{\ell, x}  = \mrm{T}_{\mrm{lb}} (\ovl{p}_x, N_{\ell, x}, \eta), \qu  d_{\ell, x} =  \mrm{T}_{\mrm{ub}} (\udl{p}_x, N_{\ell, x}, \eta),
\qu c_{\ell, y} = \mrm{T}_{\mrm{lb}} (\ovl{p}_y, N_{\ell, y}, \eta), \qu  d_{\ell, y} =  \mrm{T}_{\mrm{ub}} (\udl{p}_y, N_{\ell, y}, \eta) \eee}
for $\ell = 1, \cd, s$, where $\mrm{T}_{\mrm{lb}} (., ., .)$ and $\mrm{T}_{\mrm{ub}} (., ., .)$ are multivariate functions such that {\small
\bee &  & \mrm{T}_{\mrm{lb}} (\se, n, \eta) = \max \li \{0, \; \f{1}{n} \li \lc n \se + \f{ 1 - 2 \se - \sq{ 1 + \f{18 n \se
(1-\se)}{ \ln \f{2}{\eta}} } } {\f{2}{3n} + \f{3}{ \ln \f{2}{\eta}}}  \ri \rc  \ri \}, \\
&  & \mrm{T}_{\mrm{ub}} (\se, n, \eta) = \min \li \{1, \;  \f{1}{n} \li \lf n \se + \f{ 1 - 2 \se + \sq{ 1 + \f{18 n \se (1-\se)}{ \ln
\f{2}{\eta}} } } {\f{2}{3n} + \f{3}{ \ln \f{2}{\eta}}}  \ri \rf  \ri \} \eee} for $\se \in [0, 1], \; \eta \in (0, 1)$ and $n \in \bb{N}$. By
virtue of (\ref{vip898}), Theorem 3 of \cite{ChenTr}, and Bonferroni's inequality,  we have {\small \bel &  & \Pr \{ \tx{Rejecting}
\; \mscr{H}_i \mid \bs{p} \} \nonumber\\
 & \leq & 2 s \eta + \Pr \{ a_{\ell, x} \leq \wh{\bs{p}}_{\ell, x} \leq b_{\ell, x}, \; a_{\ell, y} \leq
\wh{\bs{p}}_{\ell, y} \leq b_{\ell, y}, \; \ell = 1, \cd, s; \; \tx{Rejecting} \; \mscr{H}_i \mid \bs{p}
\} \nonumber\\
& = & 2 s \eta + \sum_{j = 1}^s \Pr \{ a_{\ell, x} \leq \wh{\bs{p}}_{\ell, x} \leq b_{\ell, x}, \; a_{\ell, y} \leq \wh{\bs{p}}_{\ell, y} \leq
b_{\ell, y}, \; \ell = 1, \cd, s; \; \nonumber \\
&   & \qqu \qqu \qqu \wh{\bs{p}}_\ell \in \mscr{C}_\ell, \; \ell = 1, \cd, j - 1; \; \wh{\bs{p}}_j \in \mscr{D}_{i, j}^r \mid \bs{p} \} \nonumber\\
& \leq & 2 s \eta + \sum_{j = 1}^s \Pr \{ a_{\ell, x} \leq \wh{\bs{p}}_{\ell, x} \leq b_{\ell, x}, \; a_{\ell, y} \leq \wh{\bs{p}}_{\ell, y}
\leq b_{\ell, y}, \; \ell = 1, \cd, j; \; \nonumber \\
&   & \qqu  \qqu \qqu \wh{\bs{p}}_\ell \in \mscr{C}_\ell, \; \ell = 1, \cd, j - 1; \; \wh{\bs{p}}_j \in \mscr{D}_{i, j}^r \mid \bs{p} \}
\la{hi8} \eel} and {\small \bel \Pr \{ \tx{Rejecting} \; \mscr{H}_i \mid \bs{p} \} & \geq & \sum_{j = 1}^s \Pr \{ c_{\ell, x} \leq
\wh{\bs{p}}_{\ell, x} \leq d_{\ell, x}, \; c_{\ell, y} \leq \wh{\bs{p}}_{\ell, y}
\leq d_{\ell, y}, \; \ell = 1, \cd, j; \;  \nonumber \\
&   & \qqu \qu \wh{\bs{p}}_\ell \in \mscr{C}_\ell, \; \ell = 1, \cd, j - 1; \; \wh{\bs{p}}_j \in \mscr{D}_{i, j}^r \mid \bs{p}  \} \la{hi9}
\eel} for all $\bs{p} \in \mcal{Q}$ and $i = 0, 1, \cd, m-1$.  It can be seen that the bounds of $\Pr \{ \tx{Rejecting} \; \mscr{H}_i \mid
\bs{p} \}$ in (\ref{hi8}) and (\ref{hi9}) can be expressed as summations of probabilistic terms like $\Pr \{  \bs{K}_i \in \mscr{K}_i, \; 1 \leq
i \leq \ell \mid \bs{p} \}, \; \ell = 1, \cd, s$, where $\mscr{K}_\ell$ does not depend on $\bs{p}$ and is a subset of the support of
$\bs{K}_\ell$ for $\ell = 1, \cd, s$.  Since $\mscr{K}_\ell$ is independent of $\bs{p}$, we can employ the recursive technique for bounding $\Pr
\{  \bs{K}_i \in \mscr{K}_i, \; 1 \leq i \leq \ell \mid \bs{p} \}, \; \ell = 1, \cd, s$ over $\mcal{Q}$, which has been established in
\cite{Chen_EST}.

Given that the Adapted Branch and Bound algorithm can be used to determine whether $\Pr \{ \tx{Rejecting} \; \mscr{H}_i \mid \bs{p} \}$ is no
greater $\de_i$ for $\bs{p} \in \{ (p_x, p_y): 0 < p_x < 1, \; 0 < p_y < 1, \; g(p_x, p_y) \in \varTheta_i \}$,  the bisection risk tuning
technique proposed in \cite{Chen_EST, Chen_TH} can be applied to determine $\ze$ as large as possible such that the risk requirement
(\ref{mainreq}) can be satisfied.

\section{Conclusion}

In this paper, we have developed a general approach for testing composite parametric hypotheses.  We demonstrate that the problem of testing
multiple hypotheses can be cast into a unified framework of constructing sequential random intervals with prescribed coverage probabilities. We
propose inclusion principle for constructing multistage testing plans.  Instead of using likelihood ratio to define testing plans, we use
confidence limits to make early stopping possible.  For a wide class problems of testing population means, the sample sizes of our proposed
testing plans are absolutely bounded. Moreover, our testing plans are more efficient over a wide range of parametric values as compared to SPRT
and its variations. Our approach can be applied to test arbitrary number of composite hypotheses.

\appendix

\section{Preliminary Results} \la{pre}

We need some preliminary results.

\beL

\la{ULE_Basic}

Let $\mscr{E}$ be an event determined by random tuple $(X_1, \cd,
X_{\mbf{r}})$. Let $\varphi(X_1, \cd, X_{\mbf{r}})$ be a ULE of
$\se$. Then,

(i) $\Pr \{ \mscr{E} \mid \se \}$ is non-increasing with respect to
$\se \in {\Se}$ no less than $z$ provided that $\mscr{E} \subseteq
\{ \varphi(X_1, \cd, X_{\mbf{r}}) \leq z \}$.

(ii) $\Pr \{ \mscr{E} \mid \se \}$ is non-decreasing with respect to
$\se \in {\Se}$ no greater than $z$ provided that $\mscr{E}
\subseteq \{ \varphi(X_1, \cd, X_{\mbf{r}}) \geq z \}$.

\eeL

\bpf

We first consider the case that $X_1, X_2, \cd$ are discrete random
variables. Let $I_{\mbf{r}}$ denote the support of $\mbf{r}$, i.e.,
$I_{\mbf{r}} = \{ \mbf{r} (\om): \om \in \Om \}$. Define $\mscr{X}_r
= \{ (X_1 (\om), \cd, X_r(\om) ): \om \in \mscr{E}, \; \mbf{r} (\om)
= r \}$ for $r \in I_{\mbf{r}}$. Then, \be \la{eeewhole}
 \Pr \{ \mscr{E} \mid \se \}  = \sum_{r \in I_{\mbf{r}}} \;
 \sum_{ (x_1, \cd, x_r) \in \mscr{X}_r } \Pr \{ X_i = x_i, \; i =
1, \cd, r \mid \se \}. \ee

To show statement (i), using the assumption that $\mscr{E} \subseteq
\{ \varphi(X_1, \cd, X_{\mbf{r}}) \leq z \}$, we have $\varphi(x_1,
\cd, x_r) \leq z$ for $(x_1, \cd, x_r) \in \mscr{X}_r$ with $r \in
I_{\mbf{r}}$.  Since $\varphi(X_1, \cd, X_{\mbf{r}})$ is a ULE of
$\se$, we have that $\Pr \{ X_i = x_i, \; i = 1, \cd, r \mid \se \}$
is non-increasing with respect to $\se \in \Se$ no less than $z$. It
follows immediately from (\ref{eeewhole}) that statement (i) is
true.

To show statement (ii), using the assumption that $\mscr{E}
\subseteq \{ \varphi(X_1, \cd, X_{\mbf{r}}) \geq z \}$, we have
$\varphi(x_1, \cd, x_r) \geq z$ for $(x_1, \cd, x_r) \in \mscr{X}_r$
with $r \in I_{\mbf{r}}$.  Since $\varphi(X_1, \cd, X_{\mbf{r}})$ is
a ULE of $\se$, we have that $\Pr \{ X_i = x_i, \; i = 1, \cd, r
\mid \se \}$ is non-decreasing with respect to $\se \in \Se$ no
greater than $z$. It follows immediately from (\ref{eeewhole}) that
statement (ii) is true.

For the case that $X_1, X_2, \cd$ are continuous random variables,
we can also show the lemma by modifying the argument for the
discrete case. Specially, the summation of likelihood function $\Pr
\{ X_i = x_i, \; i = 1, \cd, r \mid \se \}$ over the set of tuple
$(x_1, \cd, x_r )$ is replaced by the integration of the joint
probability density function $f_{X_1, \cd, X_r} (x_1, \cd, x_r,
\se)$ over the set of $(x_1, \cd, x_r)$. This concludes the proof of
Lemma \ref{ULE_Basic}.

\epf

 \beL \la{ProbTrans}  $\Pr \{ F_Z(Z) \leq \al \}
\leq \al$ and $\Pr \{ G_Z(Z) \leq \al \} \leq \al$ for any random
variable $Z$ and positive number $\al$, where $F_Z(z) = \Pr \{  Z
\leq z \}$ and $G_Z(z) = \Pr \{ Z \geq z \}$. \eeL

Actually, Lemma \ref{ProbTrans} is a well-known result, which can be proved by the following argument.  Let $I_Z$ denote the support of $Z$.  If
$\{ z \in I_Z: F_Z(z) \leq \al \}$ is empty, then, $\{ F_Z(Z) \leq \al \}$ is an impossible event and thus $\Pr \{ F_Z(Z) \leq \al \} = 0$.
Otherwise, we can define $z^\star = \max \{ z \in I_Z: F_Z(z) \leq \al \}$.  It follows from the definition of $z^\star$ that $F_Z (z^\star )
\leq \al$. Since $F_Z(z)$ is non-decreasing with respect to $z$, we have $\{ F_Z(Z) \leq \al \} = \{ Z \leq z^\star \}$. Therefore, $\Pr \{
F_Z(Z) \leq \al \} = \Pr \{ Z \leq z^\star \} = F_Z (z^\star ) \leq \al $ for any $\al > 0$. By a similar method, one can show $\Pr \{ G_Z(Z)
\leq \al \} \leq \al$ for any $\al > 0$.

\beL \la{Unify_CH} Let $X$ be a random variable parameterized by its
mean $\bb{E} [ X ] = \se \in \Se$.  Suppose that $X$ is a ULE of
$\se$. Let $\ovl{X}_n = \f{ \sum_{i = 1}^n X_i } {n}$, where $X_1,
\cd, X_n$ are  i.i.d. samples of random variable $X$.  Then, \bee &
& \Pr
\{ \ovl{X}_n \leq z \} \leq [ \mscr{C} (z, \se) ]^n, \qqu \fa z \leq \se\\
&  & \Pr \{ \ovl{X}_n \geq z \} \leq [ \mscr{C} (z, \se) ]^n, \qqu
\fa z \geq \se. \eee Moreover, $\mscr{C} (z, \se)$ is non-decreasing
with respect to $\se$ no greater than $z$ and is non-increasing with
respect to $\se$ no less than $z$.  Similarly,  $\mscr{C} (z, \se)$
is non-decreasing with respect to $z$ no greater than $\se$ and is
non-increasing with respect to $z$ no less than $\se$.

\eeL

\bpf By the convexity of function $e^x$ and Jensen's inequality, we
have $\inf_{\ro > 0 } \bb{E} [e^{\ro (X - z)}] \geq \inf_{\ro > 0 }
e^{\ro \bb{E} [ X - z]} \geq 1$ for $\se \geq z$. In view of
$\inf_{\ro \leq 0 } \bb{E} [e^{\ro (X - z)}] \leq 1$, we have
$\mscr{C} (z, \se) = \inf_{\ro \leq 0 } \bb{E} [e^{\ro (X - z)}]$
for $\se \geq z$. Clearly, $\mscr{C} (z, \se) = \inf_{\ro \leq 0 }
e^{- \ro z} \bb{E} [e^{\ro X}]$ is non-decreasing with respect to
$z$ less than $\se$. Since $X$ is a ULE of $\se$, we have that
$\bb{E} [e^{\ro (X - z)}] = e^{- \ro z} \bb{E} [e^{\ro X} ] = e^{-
\ro z} \int_{u = 0}^\iy \Pr \{ e^{\ro X}
> u \} du$ is non-increasing with respect to $\se \geq z$ for $\ro
\leq 0$ and thus $\mscr{C} (z, \se)$ is non-increasing with respect
to $\se$ greater than $z$.

Observing that $\inf_{\ro \geq 0 } \bb{E} [e^{\ro (X - z)}] \leq 1$
and that $\inf_{\ro < 0 } \bb{E} [e^{\ro (X - z)}] \geq \inf_{\ro <
0 } e^{\ro \bb{E} [ X - z]} \geq 1$ for $\se < z$, we have $\mscr{C}
(z, \se) = \inf_{\ro \geq 0 } \bb{E} [e^{\ro (X - z)}]$ for $\se <
z$. Clearly, $\mscr{C} (z, \se) = \inf_{\ro \geq 0 } e^{- \ro z}
\bb{E} [e^{\ro X}]$ is non-increasing with respect to $z$ greater
than $\se$. Since $X$ is a ULE of $\se$, we have that $\bb{E}
[e^{\ro (X - z)}] = e^{- \ro z} \int_{u = 0}^\iy \Pr \{ e^{\ro X} >
u \} du$ is non-decreasing with respect to $\se$ for $\ro > 0$ and
consequently,  $\mscr{C} (z, \se)$ is non-decreasing with respect to
$\se$ smaller than $z$.

Making use of the established fact $\inf_{\ro \leq 0 } \bb{E}
[e^{\ro (X - z)}] = \mscr{C} (z, \se)$ and the Chernoff bound $\Pr
\{ \ovl{X}_n \leq z \}  \leq \li [ \inf_{\ro \leq 0 } \bb{E} [e^{\ro
(X - z)}] \ri ]^n$ (see, \cite{Chernoff}), we have $\Pr \{ \ovl{X}_n
\leq z \} \leq \li [ \mscr{C} (z, \se) \ri ]^n$ for $z \leq \se$.
Making use of the established fact $\inf_{\ro \geq 0 } \bb{E}
[e^{\ro (X - z)}] = \mscr{C} (z, \se)$ and the Chernoff bound $\Pr
\{ \ovl{X}_n \geq z \} \leq \li [ \inf_{\ro \geq 0 } \bb{E} [e^{\ro
(X - z)}] \ri ]^n$, we have $\Pr \{ \ovl{X}_n \geq z \} \leq [
\mscr{C} (z, \se) ]^n$ for $z \geq \se$. This concludes the proof of
Lemma \ref{Unify_CH}. \epf

\section{Proof of Theorem 1 }  \la{thm1_app}

Since $\ovl{X}_n$ is a ULE for $\se$, it follows from Lemma
\ref{ULE_Basic} that $G_{\ovl{X}_n} (z, \se)$ is non-decreasing with
respect to $\se \in \Se$. Hence, by the definition of the lower
confidence limit, $L( \f{\se_0 + \se_1}{2}, \ze \al ) \geq \se_0$ if
$G_{\ovl{X}_n} ( \f{\se_0 + \se_1}{2}, \se_0 ) \leq \ze \al$.  Since
$\ovl{X}_n$ is a ULE for $\se$, it follows from Lemma \ref{Unify_CH}
that
\[
G_{\ovl{X}_n} \li ( \f{\se_0 + \se_1}{2}, \se_0 \ri ) = \Pr \li \{
\ovl{X}_n \geq \f{ \se_1 + \se_0 }{2}  \mid \se_0 \ri \} \leq \li [
\mscr{C} \li ( \f{\se_0 + \se_1}{2}, \se_0 \ri ) \ri ]^n \leq \ze
\al
\]
if \be \la{cona} n \geq \f{ \ln (\ze \al) } { \ln \mscr{C} (
\f{\se_0 + \se_1}{2}, \se_0 ) }. \ee Observing that $G_{\ovl{X}_n}
(z, \se)$ is non-increasing with respect to $z$, we have that  $L(z,
\ze \al)$ is non-decreasing with respect to $z$.  This implies that
$L(z, \ze \al) \geq L( \f{\se_0 + \se_1}{2}, \ze \al ) \geq \se_0$
for $z \geq \f{\se_0 + \se_1}{2}$ if (\ref{cona}) holds.

On the other hand, since $\ovl{X}_n$ is a ULE for $\se$, it follows
from Lemma \ref{ULE_Basic} that $F_{\ovl{X}_n} (z, \se)$ is
non-increasing with respect to $\se$. Hence, by the definition of
the upper confidence limit, $U( \f{\se_0 + \se_1}{2}, \ze \ba ) \leq
\se_1$ if $F_{\ovl{X}_n} ( \f{\se_0 + \se_1}{2}, \se_1 ) \leq \ze
\ba$. Since $\ovl{X}_n$ is a ULE for $\se$, it follows from Lemma
\ref{Unify_CH} that
\[
F_{\ovl{X}_n} \li ( \f{\se_0 + \se_1}{2}, \se_1 \ri ) = \Pr \li \{
\ovl{X}_n \leq \f{ \se_1 + \se_0 }{2}  \mid \se_1 \ri \} \leq \li [
\mscr{C} \li ( \f{\se_0 + \se_1}{2}, \se_1 \ri ) \ri ]^n \leq \ze
\ba
\]
if \be \la{conb}
 n \geq \f{ \ln (\ze \ba) } {\ln \mscr{C} ( \f{\se_0 +
\se_1}{2}, \se_1 ) }. \ee Observing that $F_{\ovl{X}_n} (z, \se)$ is
non-decreasing with respect to $z$, we have that  $U(z, \ze \al)$ is
non-decreasing with respect to $z$.   This implies that $U(z, \ze
\ba) \leq U( \f{\se_0 + \se_1}{2}, \ze \ba ) \leq \se_1$ for $z \leq
\f{\se_0 + \se_1}{2}$ if (\ref{conb}) holds.

Let $n_{\mrm{max}}$ be the smallest integer $n$ such that
 $\{ L(\ovl{X}_n, \ze \al) < \se_0 < \se_1 < U(\ovl{X}_n, \ze \ba) \} = \emptyset$.
From the above analysis, we know that either $U(z, \ze \ba) \leq
\se_1$ or $L(z, \ze \al) \geq \se_0$ must be true for any $z \in
I_{\ovl{X}_n}$ if (\ref{cona}) and (\ref{conb}) are satisfied. This
implies that \be \la{crit} \mbf{n} \leq n_{\mrm{max}} < \max \li \{
\f{ \ln (\ze \al) } { \ln \mscr{C} ( \f{\se_0 + \se_1}{2}, \se_0 )
}, \; \f{ \ln (\ze \ba) } { \ln \mscr{C} ( \f{\se_0 + \se_1}{2},
\se_1 ) } \ri \}. \ee By the definition of the testing plan,  we
have \bee
 \Pr \{ \tx{Reject}
\; \mscr{H}_0 \mid \se_0 \} & \leq & \sum_{n = 1}^{ n_{\mrm{max}}}
\Pr \{ L (\ovl{X}_n , \ze \al) \geq \se_0, \; \mbf{n} = n \mid \se_0
\}\\
& \leq & \sum_{n = 1}^{ n_{\mrm{max}}} \Pr \{ L (\ovl{X}_n , \ze
\al) \geq \se_0 \mid \se_0
\}\\
& \leq & n_{\mrm{max}} \; \ze \al \eee and \bee
 \Pr \{ \tx{Accept}
\; \mscr{H}_0 \mid \se_1 \} & \leq & \sum_{n = 1}^{ n_{\mrm{max}}}
\Pr \{ U (\ovl{X}_n , \ze \ba) \leq \se_1, \; \mbf{n} = n \mid \se_1
\}\\
& \leq & \sum_{n = 1}^{ n_{\mrm{max}}} \Pr \{ U (\ovl{X}_n , \ze
\ba) \leq \se_1 \mid \se_1
\}\\
& \leq & n_{\mrm{max}} \; \ze \ba. \eee  Making use of (\ref{crit}),
we have that both $n_{\mrm{max}} \ze \al$ and $n_{\mrm{max}} \ze
\ba$ tend to $0$ as $\ze \to 0$.  It follows that $\Pr \{
\tx{Reject} \; \mscr{H}_0 \mid \se_0 \}$ and $\Pr \{ \tx{Accept} \;
\mscr{H}_0 \mid \se_1 \}$ tend to $0$ as $\ze \to 0$.  From the
definition of the testing plan, we have that $\{ \tx{Reject} \;
\mscr{H}_0 \}$ and $\{ \tx{Accept} \; \mscr{H}_0 \}$ are determined
by random tuple $(X_1, \cd, X_{\mbf{n}})$.  Moreover, $\{
\tx{Reject} \; \mscr{H}_0 \} \subseteq \{  \ovl{X}_n \geq \se_0  \}$
and $\{ \tx{Accept} \; \mscr{H}_0 \} \subseteq \{  \ovl{X}_n \leq
\se_1  \}$. Hence, by Lemma \ref{ULE_Basic}, we have that $\Pr \{
\tx{Accept} \; \mscr{H}_0 \mid \se \}$ is non-increasing with
respect to $\se \in \Se$ such that $\se \notin (\se_0, \se_1)$.
This completes the proof of the theorem.

\section{Proof of Theorem 2 } \la{thm2_app}

By the definition of the lower confidence limit, $\li \{ [ \mscr{C}
( \ovl{X}_n, \mcal{L} (\ovl{X}_n, \de) ) ]^n \leq \de , \; \mcal{L}
(\ovl{X}_n, \de) \leq \ovl{X}_n \ri \}$ is a sure event. It follows
from Lemmas \ref{ProbTrans} and \ref{Unify_CH} that {\small \bee \Pr
\{ \mcal{L} (\ovl{X}_n, \de) \geq \se \}  & = & \Pr \li \{ [
\mscr{C} ( \ovl{X}_n, \mcal{L} (\ovl{X}_n, \de) ) ]^n \leq \de , \;
\se \leq \mcal{L} (\ovl{X}_n, \de) \leq \ovl{X}_n \ri \}\\
& \leq & \Pr \li \{ G_{\ovl{X}_n} ( \ovl{X}_n, \mcal{L} (\ovl{X}_n,
\de) ) \leq \de, \; \se \leq \mcal{L} (\ovl{X}_n, \de) \leq
\ovl{X}_n \ri \} \leq \Pr \li \{ G_{\ovl{X}_n} ( \ovl{X}_n, \se)
\leq \de \ri \} \leq \de \eee} for $\se \in \Se$. By the definition
of the upper confidence limit, $\li \{ [ \mscr{C} ( \ovl{X}_n,
\mcal{U} (\ovl{X}_n, \de) ) ]^n \leq \de, \; \mcal{U} (\ovl{X}_n,
\de) \geq \ovl{X}_n \ri \}$ is a sure event. It follows from Lemmas
\ref{ProbTrans} and \ref{Unify_CH} that {\small \bee \Pr \{ \mcal{U}
(\ovl{X}_n, \de) \leq \se \}  & = & \Pr \li \{ [ \mscr{C} (
\ovl{X}_n, \mcal{U} (\ovl{X}_n, \de) ) ]^n
\leq \de, \; \ovl{X}_n \leq \mcal{U} (\ovl{X}_n, \de) \leq \se \ri \}\\
& \leq & \Pr \li \{ F_{\ovl{X}_n} ( \ovl{X}_n, \mcal{U} (\ovl{X}_n,
\de) )  \leq \de, \; \ovl{X}_n \leq \mcal{U} (\ovl{X}_n, \de) \leq
\se \ri \} \leq \Pr \li \{ F_{\ovl{X}_n} ( \ovl{X}_n, \se) \leq \de
\ri \} \leq \de \eee} for $\se \in \Se$.

Now, we show that $\mcal{L} (z, \de)$ is non-decreasing with respect
to $z$.  Let $z_1 < z_2$ and $L_1 = \mcal{L} (z_1, \de), \; L_2 =
\mcal{L} (z_2, \de)$.  To show the monotonicity of $\mcal{L} (z,
\de)$, it suffices to show $L_2 \geq L_1$ for three cases (not
necessarily mutually exclusive) as follows:

Case (i):  $[ \mscr{C} (z_2, z_2) ]^n \leq \de$;

Case (ii): $[ \mscr{C} (z_1, z_1) ]^n \leq \de$;

Case (iii): $[ \mscr{C} (z_1, z_1) ]^n > \de$ and $[ \mscr{C} (z_2,
z_2) ]^n > \de$.

\bsk

In Case (i),  we have $L_2 = z_2 > z_1 \geq L_1$.

In Case (ii),  we have $L_1 = z_1$ and $[ \mscr{C} (z_2, z_1) ]^n
\leq [ \mscr{C} (z_1, z_1) ]^n \leq \de$, which implies that $L_2
\geq z_1 = L_1$.

In Case (iii), it must be true that $[ \mscr{C} (z_1, L_1) ]^n =
\de$ and $[ \mscr{C} (z_2, L_2) ]^n = \de$. Suppose, to get a
contradiction,  that $L_1 > L_2$. Then, $z_1 > L_1 > L_2$ and
$\mscr{C} (z_1, L_2) < \mscr{C} (z_1, L_1) = \mscr{C} (z_2, L_2)$,
since $\mscr{C} (z, \se)$ is non-decreasing with respect to $\se
\leq z$.  On the other hand, $\mscr{C} (z_1, L_2) > \mscr{C} (z_2,
L_2)$, since $z_2 > z_1 > L_2$ and $\mscr{C} (z, \se)$ is
non-increasing with respect to $z \geq \se$.  This is a
contradiction.  Thus, it must be true that $L_1 \leq L_2$.  This
completes the proof that $\mcal{L} (z, \de)$ is non-decreasing with
respect to $z$.

In a similar manner, we can show that $\mcal{U} (z, \de)$ is
non-decreasing with respect to $z$. This completes the proof of the
theorem.

\section{Proof of Theorem 3 } \la{thm3_app}

Since $\ovl{X}_n$ is a ULE for $\se$, it follows from Lemma
\ref{Unify_CH} that $\mcal{C} (z, \se)$ is non-decreasing with
respect to $\se$ no greater than $z$. Hence, it follows from the
definition of the lower confidence limit that $\mcal{L} ( \f{\se_0 +
\se_1}{2}, \ze \al ) \geq \se_0$ if
\[
\li [ \mscr{C} \li ( \f{\se_0 + \se_1}{2}, \se_0 \ri ) \ri ]^n \leq
\ze \al, \qqu \se_0 \leq \f{\se_0 + \se_1}{2},
\]
which is guaranteed by \be \la{conc}
 n \geq \f{ \ln (\ze \al) } {
\ln \mscr{C} ( \f{\se_0 + \se_1}{2}, \se_0 ) }. \ee By Theorem
\ref{thm2}, $\mcal{L} (z, \ze \al)$ is non-decreasing with respect
to $z$. Hence, $\mcal{L} (z, \ze \al) \geq \mcal{L} ( \f{\se_0 +
\se_1}{2}, \ze \al ) \geq \se_0$ for $z \geq \f{\se_0 + \se_1}{2}$
if (\ref{conc}) holds.

On the other hand, since $\ovl{X}_n$ is a ULE for $\se$, it follows
from Lemma \ref{Unify_CH} that $\mcal{C} (z, \se)$ is non-increasing
with respect to $\se$ no less than $z$. Hence, it follows from the
definition of the upper confidence limit that $\mcal{U} ( \f{\se_0 +
\se_1}{2}, \ze \ba ) \leq \se_1$ if
\[
\li [ \mscr{C} \li ( \f{\se_0 + \se_1}{2}, \se_1 \ri ) \ri ]^n \leq
\ze \ba, \qqu \se_1 \geq \f{\se_0 + \se_1}{2},
\]
which is ensured by \be \la{cond} n \geq \f{ \ln (\ze \ba) } { \ln
\mscr{C} ( \f{\se_0 + \se_1}{2}, \se_1 ) }. \ee   By Theorem
\ref{thm2}, $\mcal{U} (z, \ze \ba)$ is non-decreasing with respect
to $z$. Consequently,
\[
\mcal{U} (z, \ze \ba) \leq \mcal{U} \li ( \f{\se_0 + \se_1}{2}, \ze
\ba \ri ) \leq \se_1 \]
 for $z \leq \f{\se_0 + \se_1}{2}$ if (\ref{cond})
holds.  This shows that either $\mcal{L} (\ovl{X}_n, \ze \al) \geq
\se_0$ or $\mcal{U} (\ovl{X}_n, \ze \ba) \leq \se_1$ must be true if
 \[
n \geq \f{ \ln (\ze \al) } { \ln \mscr{C} (  \f{\se_0 + \se_1}{2},
\se_0 ) }, \qqu  n \geq \f{ \ln (\ze \ba) } { \ln \mscr{C} (
\f{\se_0 + \se_1}{2}, \se_1 ) }.
 \]
Let $n_{\mrm{max}}$ be the smallest integer $n$ such that
 $\{ \mcal{L} (\ovl{X}_n, \ze \al) < \se_0 < \se_1 <
 \mcal{U} (\ovl{X}_n, \ze \ba) \} = \emptyset$.
 Then,  \be
 \la{cri}
  \mbf{n} \leq n_{\mrm{max}} < \max \li \{ \f{ \ln (\ze \al)
} { \ln \mscr{C} (  \f{\se_0 + \se_1}{2}, \se_0 ) }, \; \f{ \ln (\ze
\ba) } { \ln \mscr{C} (  \f{\se_0 + \se_1}{2}, \se_1 ) } \ri \}. \ee
By the definition of the test plan,  we have {\small \[
 \Pr \{ \tx{Reject}
\; \mscr{H}_0 \mid \se_0 \} \leq \sum_{n = 1}^{ n_{\mrm{max}}} \Pr
\{ \mcal{L} (\ovl{X}_n , \ze \al) \geq \se_0, \; \mbf{n} = n \mid
\se_0 \} \leq  \sum_{n = 1}^{ n_{\mrm{max}}} \Pr \{ \mcal{L}
(\ovl{X}_n , \ze \al) \geq \se_0 \mid \se_0 \} \leq n_{\mrm{max}} \;
\ze \al \]} and {\small \[
 \Pr \{ \tx{Accept}
\; \mscr{H}_0 \mid \se_1 \} \leq \sum_{n = 1}^{ n_{\mrm{max}}} \Pr
\{ \mcal{U} (\ovl{X}_n , \ze \ba) \leq \se_1, \; \mbf{n} = n \mid
\se_1 \} \leq \sum_{n = 1}^{ n_{\mrm{max}}} \Pr \{ \mcal{U}
(\ovl{X}_n , \ze \ba) \leq \se_1 \mid \se_1 \} \leq n_{\mrm{max}} \;
\ze \ba. \]}  Making use of (\ref{cri}), we have that both
$n_{\mrm{max}} \ze \al$ and $n_{\mrm{max}} \ze \ba$ tend to $0$ as
$\ze \to 0$.  It follows that $\Pr \{ \tx{Reject} \; \mscr{H}_0 \mid
\se_0 \}$ and $\Pr \{ \tx{Accept} \; \mscr{H}_0 \mid \se_1 \}$ tend
to $0$ as $\ze \to 0$.

From the definition of the testing plan, we have that $\{
\tx{Reject} \; \mscr{H}_0 \}$ and $\{ \tx{Accept} \; \mscr{H}_0 \}$
are determined by random tuple $(X_1, \cd, X_{\mbf{n}})$. Moreover,
$\{ \tx{Reject} \; \mscr{H}_0 \} \subseteq \{  \ovl{X}_n \geq \se_0
\}$ and $\{ \tx{Accept} \; \mscr{H}_0 \} \subseteq \{ \ovl{X}_n \leq
\se_1  \}$. Hence, by Lemma \ref{ULE_Basic}, we have that $\Pr \{
\tx{Accept} \; \mscr{H}_0 \mid \se \}$ is non-increasing with
respect to $\se \in \Se$ such that $\se \notin (\se_0, \se_1)$. This
completes the proof of the theorem.

\end{document}